\documentstyle[12pt]{article}
\begin{document}
\parindent 16pt
\title{\Large\bf   Ekeland variational principles with set-valued objective functions and set-valued perturbations$^1$}
\setcounter{footnote}{1}
 \footnotetext{ This work was supported by
the National Natural Science Foundation of China (Grant Nos. 11471236, 11561049).
}

\author{  {QIU JingHui
}\\
{\footnotesize\sl School of Mathematical Sciences,  Soochow
University, Suzhou 215006,  China}\\
{\footnotesize\sl Email: qjhsd@sina.com}
}
\date{}
\maketitle
\vskip 0.5cm {{\bf Abstract}  \   In the setting of real vector spaces,   we establish a general  set-valued Ekeland variational principle (briefly, denoted by EVP),  where the objective function is a set-valued map taking values in a real vector space  quasi-ordered by a convex cone $K$ and the perturbation consists of a $K$-convex subset  $H$ of the ordering cone $K$  multiplied by the distance function. Here, the assumption on lower boundedness  of the objective function is taken to be the weakest kind.
From the general set-valued EVP, we deduce a number of particular versions of set-valued EVP, which extend and improve the related results in the literature. In particular, we give several  EVPs for approximately efficient solutions  in set-valued optimization, where a usual assumption for $K$-boundedness (by scalarization) of the objective function's range is removed. Moreover, still under the weakest lower boundedness  condition, we present a set-valued EVP, where the objective function is a set-valued map taking values in a quasi-ordered topological vector space and the perturbation consists of a $\sigma$-convex subset of the ordering cone multiplied by the distance function.\\

 }\\


\vskip 1cm \baselineskip 18pt

\section*{ \large\bf 1  \  Introduction  }

\hspace*{\parindent}  Since the variational principle of Ekeland [12, 13] for approximate solutions of nonconvex minimization problems appeared in 1972, there have been various
generalizations and applications of the famous principle, for example, see [8, 14, 18]. Motivated by its wide applications, many authors have been interested in extending the Ekeland variational principle (briefly, denoted by EVP)
to the case with vector-valued maps or set-valued maps, for example, see [3, 6, 8, 10, 11, 15, 18, 26, 45, 50]  and the references therein. In the beginning, the most frequently exploited versions of vector EVP is as follows: the objective function $f:\, X\rightarrow Y$  is a vector-valued map taking values in a (topological) vector space $Y$ quasi-ordered by a convex cone $K$ and the perturbation is given by a nonzero vector $k_0$  of the ordering cone $K$ multiplied by the distance function $d(\cdot, \cdot)$, i.e., its form is as $d(\cdot, \cdot)\, k_0$ (disregarding a constant coefficient), for example, see [3, 8, 18, 19, 28,  33, 40, 50, 54]. Later, Bednarczuk and Zagrodny [7] proved a vector EVP, where the perturbation is given by a convex subset $H$ of the ordering cone multiplied by the distance function, i.e., its form is as $d(\cdot, \cdot)\, H$. This generalizes the case where directions of the perturbations are singletons $k_0$. More generally, Guti\'{e}rrez, Jim\'{e}nez and Novo [22] introduced a set-valued metric, which takes values in the set family consisting of all subsets of the ordering cone and satisfies the so-called triangle inequality. By using it they gave an original approach to extending the scalar-valued EVP to a vector-valued version, where the perturbation contains a set-valued metric.
They also deduced several particular versions of EVP involving approximate solutions for vector optimization problems and presented their interesting applications to optimization.
In the above EVPs given by Bednarczuk and Zagrodny [7] and given by Guti\'{e}rrez, Jim\'{e}nez and Novo [22], the objective functions are still a vector-valued (single-valued) map; and the perturbations contain a convex subset of the ordering cone and a set-valued metric taking values in the ordering cone, respectively.

On the other hand, Ha [24] introduced a strict minimizer of a set-valued map by virtue of Kuroiwa's set optimization criterion (see [32]). Using the method of cone extensions, Ha established a new version (see [24, Theorem 3.1]) of EVP for set-valued maps, which is expressed by the existence of a strict minimizer for a perturbed set-valued optimization problem. Inspired by Ha's work, Qiu [41] obtained an improvement of the Ha's result by using Gerstewitz's functionals. In the above Ha's and Qiu's versions, the perturbations are both given by $d(\cdot, \cdot)\, k_0$; and the objective functions are both a set-valued map taking values in a locally convex Hausdorff topological vector space (briefly, denoted by a locally convex space) quasi-ordered by a convex cone.

Furthermore, Liu and Ng [35], Tammer and Z$\breve{a}$linescu [51], Khanh and Quy [29], and  Flores-Ba\'{z}an,  Guti\'{e}rrez  and Novo [16] considered more general versions of EVP, where not only the  objective function is a set-valued map, but also the perturbation contains a  set-valued metric, or a convex subset of the ordering cone, or even a family of set-valued maps satisfying certain property.
In particular, Liu and Ng [35] established several set-valued EVPs, where the objective function is a set-valued map taking values in a quasi-ordered Banach space and the perturbation is as the form $\gamma\, d(\cdot, \cdot)\,H$ or $\gamma^{\prime} d(\cdot, \cdot)\, H,\, \gamma^{\prime}\in (0, \gamma)$, where $\gamma >0$ is a constant  and $H$  is a closed convex subset of the ordering cone $K$. Using the obtained EVPs, they presented some sufficient conditions ensuring the existence of error bounds for inequality systems.  Tammer and Z$\breve{a}$linescu [51] presented new minimal point theorems in product spaces and deduced the corresponding set-valued EVPs. As special cases, they derived many of  the existing EVPs and their extensions, for example, extensions of EVPs of Isac-Tammer (See [18]) and Ha's versions (See [24]). By using a lemma on a lower closed transitive reflexive relation on metric spaces, Khanh and Quy [29] got several stronger and more general versions of EVP, which extend and improve a lot of known results.
By extending Br\'{e}zis-Browder principle to partially ordered spaces, Flores-Ba\'{z}an,  Guti\'{e}rrez  and Novo [16] established a general strong minimal point existence theorem on quasi-ordered spaces and deduced several very general set-valued EVPs, where the objective function is a set-valued map and the perturbation even involves a family of set-valued  maps satisfying the so-called ``trangle inequality" property. These general set-valued EVPs include many previous EVPs and imply some new interesting results.

As we have seen, in the original version of EVP, the two requirements on the objective function are needed. One is the lower semi-continuity of the objective function; and the other is the lower boundedness of (the image of) the objective function. Concerning the lower semi-continuity assumption, ones have found that it can be replaced by a weaker one, i.e., so-called ``sequentially lower monotony". Sometimes,  it is called ``submonotone" (See [22, 42]) or ``monotonically semi-continuous" (see [7]) or ``condition (H-4)" (See [19]). The notions of  ``lower semi-continuity" and ``sequentially lower monotony" have already been extended to the case of vector-valued maps or set-valued maps; for details, see Section 2.  Concerning the lower boundedness assumption, we have a few words to say. When the objective function $f$ is a scalar-valued function, then the lower boundedness of $f$ is clear  and definite, i.e., there exists a real number $\alpha$ such that  $f(x)\geq \alpha$ for all $x\in X$, or equivalently, $\inf \{f(x):\, x\in X\}\,>\,-\infty$. However, when the objective function $f$  is a vector-valued map, or a set-valued map taking values in a (topological) vector space  $Y$ quasi-ordered by a convex cone $K$, there are various kinds of lower boundedness.  We shall discuss the details in Section 3. We shall see that there exists a kind of lower boundedness on set-valued maps, which is the weakest.

Very recently, in the setting of a real vector space not necessarily endowed with a topology, Guti\'{e}rrez,  Novo, R\'{o}denas-Pedregosa and Tanaka [23] studied the so-called nonconvex separation functional, which is generated by a nonempty set $E\subset Y$ and a nonzero point  $q\in Y$, and   is denoted by $\varphi^q_E$. They derived  the essential properties of this functional and used them for characterizing via scalarization several kinds of solutions of vector equilibrium problems whose image space is not endowed with any topology. Inspired by their work, we further
consider the so-called generalized nonconvex  separation functional, which is generated by using a  set $Q$ in place of a point  $q$ in the above functional, and it is denoted by $\varphi^Q_E$.
Being quite different with $\varphi^q_E$, \, $\varphi^Q_E$ could no longer satisfy the sub-additivity even though $E$ is a convex cone and $Q$ is a convex set. Fortunately, we find out that the sub-additivity $\varphi^Q_E(y_1 + y_2)\,\leq\, \varphi^Q_E(y_1) + \varphi^Q_E(y_2)$  still holds if $\varphi^Q_E(y_1) < 0$ and $\varphi^Q_E(y_2) < 0$.
 By using the property of generalized nonconvex separation functional and a pre-order principle in [44], we establish a general set-valued EVP, where the objective function is a set-valued map taking values in a real vector space quasi-ordered by a convex cone $K$ and the perturbation consists of a cone-convex subset $H$ of the ordering cone $K$ multiplied by the distance function $d(\cdot, \cdot)$. It deserves attention that here  the assumption on lower boundedness of the objective function $f$ is taken to be the weakest kind.

From the general set-valued EVP, we deduce  a number of particular  versions of set-valued EVP, which extend and improve many previous results, including all the above-mentioned set-valued EVPs.
In particular,  we obtain several EVPs for $\epsilon$-efficient solutions in set-valued optimization,   which extend the related results in [22, 42] from vector-valued maps to set-valued maps. Besides, the usual assumption for $K$-boundedness (by scalarization) of  the objective function's range in [22, 42] has been removed.  Moreover, still under the weakest lower boundedness assumption of the objective function, we present a set-valued EVP, where the objective function is a set-valued map taking values in a quasi-ordered topological vector space and the perturbation consists of a $\sigma$-convex  subset of the ordering cone  multiplied by the distance function. Our results  extend and improve the related results in [7, 22, 24, 29, 34, 35, 41, 42, 44, 51].

This paper is structured as follows. In Section 2, we give  some preliminaries, which include  some basic concepts on lower semi-continuity and sequentially lower monotony of set-valued maps.  Section 3  presents various kinds of  lower boundedness for set-valued maps and investigates their relationships. In Section 4, stimulated by [23], we consider generalized nonconvex separation functionals in real vector spaces and study their properties.
In Section 5, by using the generalized nonconvex separation functionals and a pre-order principle in [44], we establish a general set-valued EVP with the weakest lower boundedness assumption for the objective function.
From this, we deduce a  number of particular versions of set-valued EVP.  In Section 6,  we introduce the $(C, \epsilon)$-efficiency concept in set-valued optimization and deduce several EVPs for $(C, \epsilon)$-efficient solutions  in set-valued optimization, which extend and improve the related results in [22, 42]. Finally, in Section 7,  still under the weakest lower boundedness assumption on the objective function, we present  a set-valued EVP, where the perturbation consists of  a $\sigma$-convex subset of the ordering cone  multipied by the distance function.\\

\section*{\large\bf  2 \  Preliminaries }

\hspace*{\parindent}  Let $X$ be a nonempty set. As in [16], a binary relation $\preceq$ on $X$ is called a pre-order if it  satisfies the transitive property; a quasi-order if it satisfies the reflexive and transitive properties; a partial order if it satisfies the antisymmetric, reflexive and transitive properties.  Now, let $Y$ be a real vector space. If  $A, B\subset Y$  are nonempty and $\alpha\in {\bf R}$, then sets $A+B$ and $\alpha\,A$ are defined as follows:
$$A+B:=\,\{z\in Y:\, \exists x\in A,\ \exists y\in B \  {\rm such\ that} \  z= x+y\},$$
$$\alpha\, A:=\,\{z\in Y:\, \exists x\in A \ {\rm such\ that}\  z=\alpha\,x\}.$$
A nonempty set $K\subset Y$  is called a cone if $\alpha\, K\subset K$ for all $\alpha\geq 0$.  A cone $K$ is  called a convex cone if  $K+K \subset K$.  A convex cone $K$ can specify a quasi-order $\leq_K$ on $Y$  as follows:
$$y_1,\, y_2 \in Y, \ \  \  \    y_1\leq_K y_2\  \    \Longleftrightarrow\  \   y_1-y_2 \in -K.$$
In this case, $K$  is also called the ordering cone or positive cone.  In the following, we always assume that $K$ is nontrivial, i.e.,  $K\not=\{0\}$ and $K\not= Y$.

Next, we give several definitions concerning cone continuity of
set-valued maps (or  vector-valued maps) and discuss the relationship between them. In this section,
 we always assume that $X$ is a topological space, $Y$  is
a topological vector space (we always assume that it is Hausdorff), and $K$ is a closed convex cone in $Y$, unless other specified.\\

{\bf Definition  2.1} (See [24, 29, 44]).  \  A set-valued map $f:\,
X\rightarrow 2^Y$ is said to be $K$-lower semi-continuous (briefly, denoted by $K$-lsc) on $X$ iff for any $b\in Y$, the set $\{x\in X:\, f(x)\cap (b-K)\,\not=\,\emptyset\}$ is closed. $f$ has $K$-closed values iff  $f(x) + K$ is closed for all $x\in X$.\\

Next, we discuss a property on maps which is strictly weaker than
the lower semi-continuity. First we consider the case of
scalar-valued functions. Let $(X,d)$ be a metric space and ${\bf
R}$ be the real number space with the usual order and with the usual
topology.

A function $f:\, X \rightarrow {\bf R}$ is said to be sequentially
lower monotone (briefly, denoted by slm) iff for any sequence
$\{x_n\}\subset X$ with $x_n\rightarrow \bar{x}$ and $f(x_{n+1})\leq
f(x_n)$, we have $f(\bar{x})\leq f(x_n),\ \forall n.$
Here, we adopt the term ``sequentially lower monotone" from [25].
Chen, Cho and Yang [9]  also considered such functions and
called them lower semi-continuous from above (briefly, denoted by
lsca)  functions. By [9], we know that slm is strictly weaker than
the lower semi-continuity even for real-valued  functions. When ${\bf R}$ is replaced by a
quasi-ordered topological vector space $(Y, \leq_K)$, where $K$ is the ordering cone,
then we have the following definition on $K$-sequentially lower
monotone vector-valued maps (see [7, 19, 22, 25, 42]).

A vector-valued map $f:\, X \rightarrow (Y, \leq_K)$ is said to be
$K$-sequentially lower monotone (briefly, denoted by $K$-slm, or
slm) iff for any sequence $\{x_n\}\subset X$ with $x_n \rightarrow
\bar{x}$ and $f(x_{n+1})\leq_K f(x_n)$, we have $f(\bar{x})\leq_K
f(x_n),\ \forall n.$

In [7], a $K$-slm map is called a monotonically semi-continuous with
respect to $K$ map, in [19] it is called a map with property (H4),
 and in [22, 42] it is called a submonotone vector-valued map. The
notion of $K$-slm maps has also been extended to set-valued maps as
follows (see [29, 41]).\\

{\bf Definition  2.2.}  \ A set-valued map $f:\, X\rightarrow  2^Y$ is
said to be $K$-sequentially lower monotone (briefly, denoted by
$K$-slm, or slm) if for any sequence $\{x_n\}\subset X$ with $x_n
\rightarrow \bar{x}$ and $f(x_n)\subset f(x_{n+1})+K$, we have $
f(x_n)\subset f(\bar{x}) + K,\ \forall n.$\\

In [29], a $K$-slm set-valued map is called a weak $K$-lower semi-continuous from above (w.$K$-lsca) set-valued map.
It is easy to see (see [41]) that a $K$-lsc set-valued map is $K$-slm. But the converse is not true. \\

 \section*{\large\bf  3   \   Various kinds of lower boundedness for set-valued maps}

\hspace*{\parindent}
 For a scalar-valued function $f:\, X\rightarrow {\bf R}$, we  have
only one notion of lower boundedness: $f$ is said to be lower
bounded if there exists $\alpha\in {\bf R}$ such that $f(x)\geq
\alpha,\ \forall x\in X$, or equivalently, $\inf\{f(x): \, x\in X\}
>-\infty$. However, for a vector-valued (or  set-valued) map,
there  are various kinds of lower boundedness. In the following, we
assume that $(X,d)$ is a metric space and $(Y,\leq_K)$ is a
 topological vector space quasi-ordered by a convex cone $K$.
 The following notions of $K$-lower boundedness and quasi $K$-lower boundedness are well known and  are widely used in vector optimization and in extending EVP to the case of vector-valued or set-valued maps.\\

{\bf Definition  3.1.}   \
A set $M\subset Y$ is said to be $K$-lower bounded iff there exists $b\in Y$ such that $M\subset b+K$. Moreover, $M$ is said to be quasi $K$-lower bounded iff there exists a bounded set $B$ such that $M\subset B+K$.
A set-valued map
$f:\, X\rightarrow 2^Y$ is said to be
$K$-lower bounded on $X$  iff $f(X)$ is $K$-lower bounded.
Moreover, $f$ is said to be quasi $K$-lower bounded  iff $f(X)$ is quasi $K$-lower bounded.\\

Since a singleton is a bounded set, every $K$-lower bounded set is quasi $K$-lower bounded and hence every $K$-lower bounded set-valued map is quasi $K$-lower bounded. But the converse is not true.\\

{\bf Example 3.2.} \  Let $Y= {\bf R}^2$ be endowed with the topology generated by the usual Euclidean distance $d$. Let $ K\,=\{(\eta_1, 0)\in Y:\, \eta_1 \geq 0\}$ and $B\,=\{(0,\eta_2)\in Y:\, -1\leq \eta_2\leq 1\}$.
Then, $K$ is a closed convex cone in $Y$ and $B$ is a bounded set in $Y$. Obviously, the set  $B+K\,=\{(\eta_1, \eta_2)\in Y:\, \eta_1\geq 0,\, -1\leq\eta_2\leq 1\}$ is a quasi $K$-lower bounded set, but it is not
 $K$-lower bounded, since  for any $b=(b_1, b_2)\in Y$, $B+K$ is not contained in $b+K$. Now, let $X$ be ${\bf R}$ endowed with the usual distance and let $f:\, X\rightarrow 2^Y$ be defined as follows:
 $$ f(x)\,=\,\{(x^2, \,\eta_2):\, -1\leq\eta_2\leq 1\},\ \ \ \   x\in X.$$
Obviously, $f(X) = B+K$ is quasi $K$-lower bounded, but it is not $K$-lower bounded.\\

Sometimes,  a $K$-lower bounded map is said to be a bounded from
below map, for example, see [19, 29];  a quasi $K$-lower bounded map
is said to be a quasi-bounded from below map, for example, see [4, 5, 29, 51]. Several authors further considered lower boundedness by
scalarizations with continuous linear functionals. This needs to
assume that there exist enough continuous linear functionals, for example, see [22, 42]. Let's
recall some basic facts on topological vector spaces.

Let $Y$ be a topological vector space and let $Y^*$ be its topological dual, i.e., the
vector space consisting of all continuous linear functionals on $Y$.
 It may happen that $Y^*=\{0\}$, i.e., there is no
nontrivial continuous linear functional, even though $Y$ is Hausdorff (See [31, pp. 157-158]). However, if $Y$ is a
locally convex Hausdorff topological vector space (briefly, denoted
by a locally convex space), then $Y^*$ is large enough so that it can separates
points in $Y$, i.e., for any  two  different points $y_1\not=y_2$ in
$Y$, there exists  $y^*\in Y^*$ such that $y^*(y_1)\not= y^*(y_2)$
(For details, see [30, 31, 52]). For any $y^*\in Y^*$, we define a continuous semi-norm $p_{y^*}$ on $Y$ as follows: $p_{y^*}(y) = |y^*(y)| \ \forall y\in Y.$ The semi-norm family $\{p_{y^*}:\, y^*\in Y^*\}$ generates a locally convex Hausdorff topology on $Y$ (See, e. g., [30, 31, 52]), which  is called the weak topology on $Y$ and denoted by $\sigma(Y, Y^*)$. For any nonempty subset $F$ of $Y^*$, the semi-norm family $\{p_{y^*}:\, y^*\in F\}$ can also generate a locally convex topology (which need not be Hausdorff) on $Y$, which is denoted by $\sigma(Y, F)$.

 The positive
polar cone of $K$ is denoted by $K^+$, i.e., $K^+ = \{y^*\in Y^*:\,
y^*(y)\geq 0,\, \forall y\in
K\}$. An element in $K^+$ is called a positive continuous linear functional on $Y$. If $Y$ is a topological vector space and ${\rm int}(K)\not=\emptyset$, then $K^+\backslash\{0\}\not=\emptyset.$
 If $Y$ is a locally convex space and $0\not\in {\rm cl}(K)$, then we also have $K^+\backslash\{0\}\not=\emptyset$.\\

{\bf Definition   3.3} (See [46]).  \   If there exists $k^*\in
 K^+\backslash\{0\}$  such that $\cup\{k^*(f(x)):\, x\in X\}$ is
 bounded from below, then $f:\, X\rightarrow
 2^Y$ is said to be $k^*$-lower bounded.\\

 {\bf Proposition    3.4} (See [46]). \  {\sl Let $f:\, X\rightarrow
 2^Y$ be  quasi
 $K$-lower bounded. Then, for any $k^*\in
 K^+\backslash\{0\}$, $f$ is $k^*$-lower bounded.}\\

Let $Y$ be a locally convex space, $K\subset Y$ be a convex cone and $H\subset Y$ be a $K$-convex set, i.e., $H+K$ is convex.
By the Hahn-Banach separation theorem, we can show that $0\not\in {\rm cl}(H+K)$ iff $K^+\cap H^{+s} \not=\emptyset$, where  $H^{+s}$ denotes the set $\{y^*\in Y^*:\,  \exists\ \delta >0 \  {\rm such \ that } \   y^*(h)\geq \delta\   \forall h\in H\}$.\\

{\bf Definition   3.5.} \  let $Y$ be a locally convex space, $K\subset Y$ be a convex cone, $H\subset K\backslash -K$ be a $K$-convex set and $0\not\in {\rm cl}(H+K)$. If  there exists $k^*\in
 K^+\cap H^{+s}$  such that $\cup\{k^*(f(x)):\, x\in X\}$ is
 bounded from below, then $f:\, X\rightarrow
 2^Y$ is said to be $k^*(H)$-lower bounded.\\

Particularly, in Definition 3.5, the set $H$ may be a singleton $\{k_0\}$, where $k_0\in K\backslash -K$ such that $k_0\not\in {\rm cl}(K)$.
The following example shows that there being $k^*\in K^+\cap H^{+s}$ such that $f$ is $k^*$-lower bounded doesn't imply that $f$ is quasi $K$-lower bounded.\\

{\bf  Example    3.6.} \  Let $X$ be ${\bf R}$ endowed with
the usual metric, $Y$ be ${\bf R}^2$ endowed with the
topology generated by the Euclidean distance, $K\subset Y$ be the convex cone $\{(y_1, y_2)\in {\bf R}^2:\, y_1\geq 0\
{\rm and}\ y_2\geq 0\}$, and $H$ be the set $\{(y_1, y_2)\in {\bf R^2}:\,y_1\geq 0,\, y_2\geq 0 \  {\rm and} \ y_1+y_2 = 1\}$.
Obviously, $H\subset K\backslash -K$ and $0\not\in {\rm cl}(H+K)$.
 Let $k^*\in Y^*$ be defined as follows: $k^*(y_1,
y_2) \,=\, y_1 + y_2$, \ $(y_1, y_2)\in Y={\bf R}^2$.
Clearly,
$$k^*(y)\,=\, y_1 +y_2 \,\geq\,0\  \  \  \forall y=(y_1, y_2)\in K;$$
and
$$k^*(y)\,=\, y_1+y_2 \,=\,1\ \  \   \forall y=(y_1, y_2)\in H.$$
Thus, $k^*\in K^+\cap H^{+s}$.
 Let $f:\, X\rightarrow
2^Y$ be defined as follows: $f(x) = \{(|x|,
x)\}\subset Y ={\bf R}^2,\ \  \   x\in X={\bf R}$. It is easy to see that
$k^*(f(x)) = |x|+x  \geq 0,\ \forall x\in X$. Thus, $f(X)$ is
$k^*$-lower bounded. But, there isn't a bounded set $B$ such that
$f(X) \subset B+K$, since $p_2\circ f(X)= (-\infty, +\infty)$, where $p_2$ is the projection: $(y_1, y_2)\mapsto y_2$.\\

From Proposition 3.4 and Example 3.6, we know that $k^*$- or
$k^*(H)$-lower boundedness is strictly weaker than the quasi
$K$-lower boundedness.
The following kind of lower boundedness seems to be the weakest; also refer to  [41, 46].\\

 {\bf Definition   3.7.} \  Let $H\subset  K\backslash -K$. A
 set-valued map $f:\, X\rightarrow 2^Y$ is
 said to be $H$-lower bounded if there exists $y_0\in Y$ and $\epsilon >0$ such
 that
$$f(X)\cap(y_0-\epsilon\,H-K)\,=\,\emptyset.$$

 {\bf Proposition  3.8} (Refer to [41, 46]). \  {\sl Let $0\not\in {\rm cl}(H+K)$. If $f:\,
 X\rightarrow 2^Y$ is quasi $K$-lower bounded,
 then for any $y\in Y$, there exists $\epsilon >0$ such that
 $$f(X)\cap(y-\epsilon\,H-K)\,=\,\emptyset.$$
  Certainly, $f$ is $H$-lower bounded.}\\

{\sl  Proof.} \ Assume the contrary. There exists $y_0\in
Y$ such that
$$f(X)\cap (y_0-n\, H -K)\,\not=\,\emptyset,\ \  \forall
n.\eqno{(3.1)}$$ Since $f$ is quasi $K$-lower bounded, there exists
a bounded set $B$ in $Y$ such that $f(X) \subset K+B$. Combining
this with (3.1), we have
$$(K+B)\cap(y_0-n\, H -K)\,\not=\,\emptyset,\ \forall n.$$
For each $n$, there exists $k_n\in K,\ b_n\in B$
 such that
$$k_n + b_n\,\in \, y_0 -n\, H - K.$$
From this,
$$\frac{b_n}{n} -\frac{y_0}{n}\,\in\, -H-K- \frac{k_n}{n}\,\subset\, -H-K.$$
 Letting $n\rightarrow\infty$, we have
$0\,\in\, -{\rm cl}(H+K)$. This contradicts $0\not\in{\rm cl}(H+K)$. \hfill\framebox[2mm]{}\\

 {\bf Proposition  3.9.}  \ {\sl Let $K^+\cap H^{+s}\,\not=\,\emptyset$ and $f:\, X\rightarrow 2^Y$ be $k^*(H)$-lower bounded. Then, for any $y\in Y$, there  exists $\epsilon >0$ such that
$$f(X) \cap (y-\epsilon\,H-K)\,=\,\emptyset.$$}

{\sl  Proof.} \  Let $k^*\in K^+\cap H^{+s}$ such that $k^*\circ f(X)$ is lower bounded.
Assume the contrary. There exists $y_0\in Y$ such that
$$f(X) \cap (y_0-n\,H-K)\,\not=\,\emptyset,\ \ \ \forall n.$$
Thus, for each $n$, there exists $y_n\in f(X)$ such that
$$y_n\,\in\, y_0-n\,H-K.\eqno{(3.2)}$$
Since $k^*\in H^{+s}\cap K^+$, there  exists $\delta >0$ such that $k^*(h)\geq \delta \ \  \forall h\in H$ and $k^*(y)\geq 0\ \  \forall y\in K.$
Combining this with (3.2), we have
$$k^*(y_n)\,\leq\, k^*(y_0) - n\,\delta,\ \ \ \forall n.$$
Here, every $y_n\in f(X)$. This is contradicts that $k^*\circ f(X)$ is lower bounded. \hfill\framebox[2mm]{}\\

The following example shows that even though for any $y\in Y$,  there exists $\epsilon >0$ such that
$$f(X) \cap (y-\epsilon\,H-K)\,=\,\emptyset,$$
it still may happen that for every $k^*\in K^+\backslash\{0\}$, $f$ is not $k^*$-lower bounded.\\

{\bf Example 3.10.} \  Let $X={\bf R}$ and let $Y={\bf R^2}$ be  endowed with the topology generated by the Euclidean distance.  Let the convex cone $K$ be $\{(y_1, y_2)\in Y:\, y_1\geq 0,\, y_2\geq 0\}$ and let $H$ be $\{(1-\lambda)(1,1) \,+\,\lambda(2,1)\,\in \,Y:\, 0\leq\lambda\leq 1\}$. It is easy to verify that
$$H+K\,=\,\{(1+\lambda +y_1,\, 1+y_2)\,\in \,Y:\, y_1\geq 0,\, y_2\geq 0,\, 0\leq \lambda \leq 1\}.$$
Clearly, $H\subset K\backslash-K$ and $0\not\in {\rm cl}(H+K)$. Define a set-valued map $f:\, X\rightarrow 2^Y$ as follows:
$$f(x)\,=\,\{(x,0),\, (0,x)\}\,\subset\,Y,\ \ \ x\in X.$$
Obviously,
$$f(X)\,=\, \{(x,0)\in Y:\, x\in {\bf R}\}\,\cup\,\{(0,x)\in Y:\,x\in {\bf R}\}.$$
For any given $y_0=(y_{01},  y_{02})\in Y$, take $\epsilon > |y_{01}| \,+\,|y_{02}|.$ Then,
\begin{eqnarray*}
y_0-\epsilon\,H-K&=&(y_{01}, y_{02})+\{(-\epsilon(1+\lambda)-y_1,\, -\epsilon -y_2):\, y_1\geq 0,\, y _2\geq 0,\, 0\leq\lambda\leq 1\}\\
&=&\{(y_{01}-\epsilon(1+\lambda)-y_1,\, y_{02}-\epsilon-y_2):\, y_1\geq 0,\, y_2\geq 0,\, 0\leq\lambda\leq1\}.
\end{eqnarray*}
Assume that $y_{01}-\epsilon(1+\lambda)-y_1\,=\,0$. Then
$$y_{01}\,=\,\epsilon(1+\lambda) + y_1\,\geq\,\epsilon(1+\lambda) \geq \epsilon,$$
contradicting  $\epsilon > |y_{01}| \geq y_{01}$. Similarly,
assume that $y_{02}-\epsilon -y_2 = 0.$ Then $y_{02} = \epsilon +y_2 \geq \epsilon,$
contradicting $\epsilon > |y_{02}| \geq y_{02}.$
Thus, $$f(X)\,=\,\{(x,0):\, x\in {\bf R}\}\,\cup\,\{(0, x):\, x\in {\bf R}\}$$ does not intersect  $y_0 -\epsilon\,H -K$.
That is,
$$f(X)\,\cap\, (y_0 -\epsilon\,H-K)\,=\,\emptyset.$$
Hence $f$ is $H$-lower bounded.
 For any $y^*\in
Y^*\backslash\{0\}$, there exists a unique $(\alpha, \beta)\in {\bf
R}^2$ such that $y^*(y_1, y_2)\,=\, \alpha \,y_1 +\beta\, y_2,\ \
\forall (y_1, y_2)\in Y={\bf R}^2$, where $|\alpha| +|\beta| >0$.
Obviously, $y^*\circ f(x) = \{\alpha\,x\}\cup\{\beta\,x\}$. Hence,
$y^*\circ f(X) = \{\alpha\,x:\, x\in {\bf R}\}\cup\{\beta\,x:\,
x\in{\bf R}\}$ is not lower bounded in real number space ${\bf R}$.
That is, $f(X)$ is not $y^*$-lower bounded.  Certainly, for any
$k^*\in K^+\backslash\{0\}$, $f$ is not $k^*$-lower bounded.\\

Summing up the main points of the section, we have the following scheme:
$$\begin{array}{clcrc}
\mbox{$K$-lower boundedness} \
\raisebox{-1ex}{$\stackrel{\Rightarrow}{\not\Leftarrow}$}\
 \mbox{quasi $K$-lower boundedness} \
 \raisebox{-1ex}{$\stackrel{\Rightarrow}{\not\Leftarrow}$}\
 \mbox{$k^*(H)$-lower boundedness}
\end{array} $$
$$\begin{array}{clcrc}
\mbox{$k^*(H)$-lower boundeness} \
\raisebox{-1ex}{$\stackrel{\Rightarrow}{\not\Leftarrow}$}\
 \mbox{$H$-lower boundedness}.
\end{array}$$

\section*{\large\bf  4   \   Generalized nonconvex separation functional}

\hspace*{\parindent}
 Very recently, Guti\'{e}rrez, Novo, R\'{o}denas-Pedregosa and Tanaka [23] studied the so-called noncnvex separation functional in a real vector space not necessarily endowed  with a topology. They derived the essential properties of this functional  and successfully applied them for characterizing via scalarization several kinds of solutions of vector equilibrium problems whose image space is not endowed with any particular topology.
As in [23], let $Y$ be a real vector space,  $q\in Y\backslash\{0\}$ and $\emptyset\not= E\subset Y$. The so-called nonconvex separation functional $\varphi_E^q:\, Y\rightarrow {\bf R}\cup\{\pm\infty\}$  is defined as follows:
$$
\varphi^q_E(y):=\left\{
\begin{array}{cc}

+\infty    &{\rm if}  \  y\not\in {\bf R}q-E,\\

\inf\{t\in {\bf R}:\, y\in t q-E\}  \   &{\rm otherwise}.

\end{array}
\right.
$$

  This functional was   introduced in [17] and it is called  by different names: Gerstewitz's function, nonlinear scalarization function, smallest strictly monotonic function [36], shortage function [37], and so on.  For its main properties, please see, for example,  [8, 18, 23, 36, 47, 48, 50] and the references therein. We observe that in [23] the definition of nonconvex separation functional is stated in  the most general form, where $E$ is an arbitrary nonempty set and $q$ is an arbitrary direction, not assuming any hypothesis on $E$  and $q$. Inspired by  Guti\'{e}rrez, Novo, R\'{o}denas-Pedregosa and Tanaka's work [23] we shall extend the nonconvex separation functional from a point $q$ into a set $Q$ and  further investigate the fundamental properties of such a generalized nonconvex separation functional.
Applying the generalized nonconvex separation functional and its properties, we shall derive EVPs  with set perturbations under the weakest  lower boundedness condition for objective functions.

First, we recall several related notions. Let $Y$ be a real vector space. For a nonempty set $A\subset Y$, the vector closure of $A$ is defined as follows (See [1, 2, 47]):
$${\rm vcl}(A)\,=\,\{y\in Y:\,  \exists v\in Y,\,\exists \lambda_n \geq 0,\,\lambda_n\rightarrow 0\ {\rm such\ that}\  y+\lambda_n v \in A,\ \forall n\in{\bf N}\}.$$
For any given $v_0\in Y$, we define the $v_0$-closure of $A$  as follows (see [43, 47]):
$${\rm vcl}_{v_0}(A) \,=\,\{y\in Y:\, \exists \lambda_n\geq  0,\, \lambda_n \rightarrow 0\ {\rm such\ that}\ y+\lambda_n v_0\in A,\ \forall n\in {\bf N}\}.$$
Obviously,
$$A\,\subset\,{\rm vcl}_{v_0}(A)\,\subset\, \cup_{v\in Y}{\rm vcl}_v(A)\,=\,{\rm vcl}(A).$$
All the above inclusions are proper (For details, see [47]).
Moreover, if $Y$ is a topological vector space and ${\rm cl}(A)$ denotes the topological closure of $A$, then ${\rm vcl}(A) \subset {\rm cl}(A)$  and the inclusion is also proper.
A subset $A$ of $Y$ is said to be $v_0$-closed iff $A = {\rm vcl}_{v_0}(A)$; to be vectorially closed iff $A = {\rm vcl}(A)$; to be (topologically) closed iff $A = {\rm cl}(A)$.
Moreover, let $\Theta\subset Y$ be a convex set.  Put $\Theta_0:=\,\cup_{0\leq\lambda\leq 1}\lambda\,\Theta$. Then $\Theta_0$ is a convex set and $0\in \Theta_0$. For $0 < \epsilon <\epsilon^{\prime}$, we have $\epsilon\,\Theta_0\,\subset\, \epsilon^{\prime}\Theta_0$. For any real sequence $\{\epsilon_n\}$ with every $\epsilon_n >0$ and $\epsilon_n\rightarrow 0$, we have   $\cap_{\epsilon >0}(A-\epsilon\,\Theta_0)\,=\,\cap_{n=1}^{\infty}(A-\epsilon_n \Theta_0).$ The set $\cap_{\epsilon >0}(A-\epsilon\,\Theta_0)$, denoted by ${\rm cl}_{\Theta}(A)$, is called the $\Theta$-closure of $A$.  $A$ is said $\Theta$-closed iff ${\rm cl}_{\Theta}(A) = A$. Particularly, if $Y$ is a locally convex space, then every locally closed set  (concerning locally closed sets, see [39]) is $\Theta$-closed for every bounded convex set $\Theta$. But, a $\Theta$-closed set, where $\Theta$ is a certain bounded convex set, may  be non-locally closed. In fact, a subset $A$ of a locally convex space $Y$ is locally closed iff for every bounded convex set $\Theta$, $A$  is $\Theta$-closed (For details, see [43]).

Inspired by [23] we introduce the so-called generalized nonconvex separation functional $\varphi^Q_E:\, Y\rightarrow {\bf R}\cup\{\pm\infty\}$ by using a set $Q$ in place of  a point $q$.\\

{\bf Definition    4.1.} \  Let $Y$ be a real vector space, $Q\subset Y\backslash\{0\}$ and $\emptyset\not= E\subset Y.$  The generalized nonconvex separation functional $\varphi^Q_E$ is defined as follows:
$$
\varphi^Q_E(y):=\left\{
\begin{array}{cc}

+\infty    &{\rm if}  \  y\not\in {\bf R}Q-E,\\

\inf\{t\in {\bf R}:\, y\in t Q-E\}  \   &{\rm otherwise}.

\end{array}
\right.
$$
If $-\infty < \varphi^Q_E(y) < +\infty$, then either $y\in \varphi^Q_E(y)\, Q  -E$ or $y\not\in \varphi^Q_E(y)\, Q-E$ and there exists  a positive sequence $(\epsilon_n)$  with $\epsilon_n \rightarrow 0$ such that $y\in (\varphi^Q_E(y) +\epsilon_n)\, Q-E.$\\

{\bf Proposition  4.2.}  \ {\sl Let $E$ satisfy $\lambda\,E \subset E$ for all $\lambda >0$ and $0\not\in  Q+E$. Then,
$$\{y\in Y:\,\varphi^Q_E(y) = -\infty\} \not= \emptyset\ \  \ \ \  \Longleftrightarrow\ \  \  \ \    0\in {\rm vcl}(Q+E).$$}

{\sl  Proof. }  \  Let $y\in Y$ such that $\varphi^Q_E(y) = -\infty$.  Then,  there exists a sequence $\{\beta_n\}$  with every $\beta_n >0$  and $\beta_n \rightarrow +\infty$  such that $y \in -\beta_n Q-E$. Thus,
$$\frac{1}{\beta_n}\, y \,\in\,-Q-E\ \  {\rm and} \ \  \frac{1}{\beta_n}\,(-y)\,\in\, Q+E.$$
From this, $0\in {\rm vcl}_{-y}(Q+E) \subset {\rm vcl}(Q+E).$

Conversely, let $0 \in {\rm vcl}(Q+E)$. Then, there exists $v\in Y$ and a sequence $\{\lambda_n\}$ with every $\lambda_n \geq 0$ and $\lambda_n \rightarrow 0$ such that $\lambda_n v \in Q+E$. Since $0\not\in Q+E$, we have every $\lambda_n >0$. Thus,
$$ v \in \frac{1}{\lambda_n} Q +E   \ \   {\rm and} \ \ -v \in -\frac{1}{\lambda_n}Q -E.$$
Put  $y=-v$. Then $\varphi^Q_E(y)  = -\infty$.\hfill\framebox[2mm]{}\\

{\bf Proposition  4.3.}  \ {\sl Let $E$ satisfy $\lambda E \subset E$ for all $\lambda \geq 0$  and $0\not\in Q+E$.

{\rm (a)} \  $\varphi^Q_E(0) \, = \,0$.

{\rm (b)} \  $\varphi^Q_E(\alpha\, y) \, = \, \alpha \,\varphi^Q_E(y) $ for all $y\in Y$ and all $\alpha \geq 0$.}\\

{\sl  Proof.} \  (a) \  Obviously, $0\in 0\cdot Q-E$, so $\varphi^Q_E(0)\leq 0$.
Assume that $\varphi^Q_E(0) < 0$.  Then, there exists $\epsilon \geq 0$ such that
$\varphi^Q_E(0) +\epsilon  <  0$ and
$$0\,\in\,(\varphi^Q_E(0) +\epsilon)Q-E\,=\,-(\varphi^Q_E(0) +\epsilon)(-Q-E).$$
From this,  $0\,\in\, -Q-E$, which contradicts the assumption tat $0\,\not\in\, Q+E$.

(b) \  We shall prove the result according to the following different cases.

Case (I) \  $\alpha =0$. By (a), $\varphi^Q_E(\alpha\,y) = \varphi^Q_E(0) = 0$.
Also, $0\cdot \varphi^Q_E(y) = 0$, where we define $0\cdot (+\infty) = 0\cdot (-\infty) = 0$ if necessary.
Hence $\varphi^Q_E(\alpha\,y)\,=\,\alpha\, \varphi^Q_E(y)$  holds when $\alpha =0$.

Case (II) \  $\alpha >0$ and $-\infty < \varphi^Q_E(y) < +\infty$. There exists a sequence $\{\epsilon_n\}$ with  every $\epsilon_n\geq 0$ and $\epsilon_n \rightarrow 0$ such that
$$y\,\in\,(\varphi^Q_E(y) + \epsilon_n) Q -E.$$
Thus,
\begin{eqnarray*}
\alpha\, y \,&\in&\,\alpha(\varphi^Q_E(y) + \epsilon_n) Q -E\\
&=&\, (\alpha\, \varphi^Q_E(y) +\alpha\,\epsilon_n) Q-E.
\end{eqnarray*}
From this,  $\varphi^Q_E(\alpha\, y)\,\leq\,\alpha\, \varphi^Q_E(y) +\alpha\,\epsilon_n,\  \forall n.$   Letting $n\rightarrow\infty$, we have
$$\varphi^Q_E(\alpha\, y)\,\leq\, \alpha\,\varphi^Q_E(y). \eqno{(4.1)}$$
On the other hand,
$$\varphi^Q_E(y) \,=\,\varphi^Q_E(\frac{1}{\alpha} \alpha \,y)\,\leq\, \frac{1}{\alpha}\,\varphi^Q_E(\alpha\,y).$$
Thus,
$$\varphi^Q_E(\alpha\,y)\,\geq\,\alpha\,\varphi^Q_E(y). \eqno{(4.2)}$$
By (4.1) and (4.2), we have
$$\varphi^Q_E(\alpha\,y)\,=\,\alpha\,\varphi^Q_E(y).$$

Case (III) \  $\alpha >0$ and $\varphi^Q_E(y) = +\infty$. For all $t\in {\bf R}$, $y\not\in t\,Q-E$. This leads to  $\alpha\,y\,\not\in\, t\,Q-E, \   \forall t\in {\bf R}.$
From this,
$$\varphi^Q_E(\alpha\, y)\,=\,+\infty\,=\,\alpha\,(+\infty)\,=\,\alpha\,\varphi^Q_E(y).$$

Case (IV) \  $\alpha>0$ and $\varphi^Q_E(y) = -\infty$. There exists a sequence $\{\beta_n\}$ with all $\beta_n >0$ and $\beta_n \rightarrow +\infty$ such that $y\,\in\, -\beta_n Q-E$. Thus,
$$\alpha\, y\,\in\, -\alpha\,\beta_n Q-E.$$
Now, $\alpha\,\beta_n \rightarrow +\infty$, we have $\varphi^Q_E(\alpha\,y)\,=\,-\infty$. Thus,
$$\varphi^Q_E(\alpha\,y)\,=\, -\infty \,=\,\alpha\,(-\infty)\,=\, \alpha\,\varphi^Q_E(y).$$
\hfill\framebox[2mm]{}\\

{\bf Proposition  4.4.}  \ {\sl Let $E$ satisfy $\lambda E\subset E$ for all $\lambda >0$.

{\rm (a)} \  Let   $Q+E$ be a convex set. Then,
$$\varphi^Q_E(y_1+y_2)\,\leq\,\varphi^Q_E(y_1)  + \varphi^Q_E(y_2)\ \ {\rm whenever} \ \  \varphi^Q_E(y_1) <0\ \ {\rm and} \ \ \varphi^Q_E(y_2) <0.$$

{\rm (b)} \   Let $Q-E$ be a convex set. Then,
$$\varphi^Q_E(y_1+y_2)\,\leq\,\varphi^Q_E(y_1)  + \varphi^Q_E(y_2)\ \ {\rm whenever} \ \  \varphi^Q_E(y_1) >0\ \ {\rm and} \ \ \varphi^Q_E(y_2) >0.$$}

{\sl  Proof.} \  (a) \  Let $-\infty < \varphi^Q_E(y_1) < 0$  and $-\infty < \varphi^Q_E(y_2)  <0$.
Then, there exists a sequence $\{\epsilon_n\}$ with all $\epsilon_n \geq 0$ and $\epsilon_n \rightarrow 0$ such that
$$  \varphi^Q_E(y_1) +\epsilon_n <0 \ \ \ \ \   {\rm and} \ \ \ \  \   y_1\,\in\,(\varphi^Q_E(y_1) + \epsilon_n)\,Q - E. \eqno{(4.3)}$$
And there exists a sequence $\{\delta_n\}$ with all $\delta_n \geq 0$ and $\delta_n \rightarrow 0$ such that
$$  \varphi^Q_E(y_2) + \delta_n <0 \ \ \ \ \  {\rm and} \ \ \  \  \  y_2\,\in\,(\varphi^Q_E(y_2) + \delta_n)\,Q-E.  \eqno{(4.4)}$$
By (4.3) and (4.4) and using that $Q+E$ is convex, we have
\begin{eqnarray*}
y_1 + y_2 \,&\in&\, (\varphi^Q_E(y_1) +\epsilon_n) Q-E +(\varphi^Q_E(y_2) + \delta_n)Q-E\\
&=&\, -(\varphi^Q_E(y_1) +\epsilon_n) (-Q-E) -(\varphi^Q_E(y_2) +\delta_n) (-Q-E)\\
&=&\, -(\varphi^Q_E(y_1) +\epsilon_n +\varphi^Q_E(y_2) + \delta_n) \,(-Q-E)\\
&\subset&\, (\varphi^Q_E(y_1) +\varphi^Q_E(y_2)  +\epsilon_n + \delta_n)\,Q - E.
\end{eqnarray*}
Hence,
$$\varphi^Q_E(y_1+y_2)\,\leq\, \varphi^Q_E(y_1)  +\varphi^Q_E(y_2)  +\epsilon_n +\delta_n$$ and   $$  \varphi^Q_E(y_1 +y_2) \,\leq\,\varphi^Q_E(y_1) + \varphi^Q_E(y_2).$$
Let one of the $\varphi^Q_E(y_1),\ \varphi^Q_E(y_2)$ be $-\infty$. For example, let $\varphi^Q_E(y_1) = -\infty$  and $\varphi^Q_E(y_2)  <0$.
By $\varphi^Q_E(y_1)  = -\infty$, there exists a sequence $\{\beta_n\}$ with all $\beta_n >0$ and $\beta_n\rightarrow\infty$ such that
$$y_1 \,\in\,-\beta_n Q-E\ \ \ {\rm and} \ \ \ y_1\,\in\, \beta_n(-Q-E).\eqno{(4.5)}$$
By $\varphi^Q_E(y_2)  < 0$, there  exists $\lambda >0$  such that
$$y_2\,\in\, -\lambda\, Q-E\,=\,\lambda(-Q-E). \eqno{(4.6)}$$
By (4.5) and (4.6) and using that $Q+E$ is convex, we have
$$y_1 +y_2 \,\in\, (\beta_n +\lambda)\,(-Q-E)\,\subset\, -(\beta_n +\lambda)\, Q-E.$$
Since $\beta_n +\lambda \rightarrow +\infty\ (n\rightarrow\infty)$, we have $\varphi^Q_E(y_1+y_2)\,=\,-\infty$.
Thus, we still have
$$\varphi^Q_E(y_1+y_2)\,\leq\, \varphi^Q_E(y_1) + \varphi^Q_E(y_2).$$

(b) \  The proof is similar to that of (a). Here we won't write the details. \hfill\framebox[2mm]{}\\

But, if $\varphi^Q_E(y_1)$ and $\varphi^Q_E(y_2)$ have different signs, for example, $\varphi^Q_E(y_1) >0$ and $\varphi^Q_E(y_2) <0$, then we don't have
$\varphi^Q_E(y_1+y_2) \,\leq\,\varphi^Q_E(y_1) + \varphi^Q_E(y_2)$ even though $E$ is a convex cone $K$  and $Q$ is a convex set $H$.\\

{\bf  Example 4.5.} \   Let $Y$ be ${\bf R}^2$ with the usual linear structure and with the ordering cone $K:=\{(\eta_1, \eta_2)\in {\bf R}^2:\, \eta_1\geq 0,\, \eta_2\geq 0\}$.
Let $H$ be the set $\{(1-\frac{\lambda}{2}, \, 1-\frac{\lambda}{2}) \in {\bf R}^2:\, 0\leq \lambda\leq 1\}$. Obviously, $H\subset K$ is a convex set and $0\not\in {\rm vcl}(H+K)$.
Put $y_1 =(1, 1)$ and put $y_2=(-1, -1)$.  It is easy to verify that
$$\varphi_K^H(y_1)\,=\,\inf\{t\in {\bf R}:\, y_1=(1, 1)\in t\, H-K\}\,=\, 1$$
and
$$\varphi_K^H(y_2)\,=\,\inf\{t\in {\bf R}:\, y_2=(-1, -1)\in t\, H-K\}\,=\, -2.$$
We remark that
$$\varphi_K^H(y_1 +y_2)  = \varphi_K^H(0, 0) =0  > 1-2=\varphi_K^H(y_1) +\varphi_K^H(y_2).$$

Now, let $E$ be a convex cone $K$ specifying a quasi-order $\leq_K$ and let $Q$ be a $K$-convex set $H$ such that $H\subset K\backslash -K$.
Obviously, $\varphi^H_K$ is nondecreasing with respect to $\leq_K$. Synthesizing the above results we have the following particular proposition, which is a convenient tool for deriving set-valued EVPs  with set-valued  perturbations (See Section 5).\\

{\bf Proposition  4.6.}  \ {\sl Let $K\subset Y$ be a convex cone and $H\subset K\backslash-K$ be a $K$-convex set. Then, the generalized nonconvex separation functional $\varphi^H_K$ has the following properties:

{\rm (a)} \  $0\,\not\in\, {\rm vcl}(H+K)\ \ \  \Longleftrightarrow\ \  \  \varphi_K^H(y)\not= -\infty , \ \forall y\in Y.$

{\rm (b)} \  $y_1\leq_K y_2\ \ \  \Longrightarrow\ \ \  \varphi_K^H(y_1)  \leq  \varphi_K^H(y_2),\  \forall y_1,\, y_2\in Y.$

{\rm (c)} \  $\varphi_K^H(\alpha\,y)\,=\,\alpha\,\varphi_K^H(y),\  \forall y\in Y,\ \forall \alpha \geq 0.$

{\rm (d)} \  $\varphi_K^H(y_1+y_2)\,\leq\,\varphi_K^H(y_1)  +\varphi_K^H(y_2)$, whenever  $\varphi_K^H(y_1)  < 0$ and  $\varphi_K^H(y_2)  < 0$.
}\\

Next, we further give two results on generalized nonconvex separation  functional with some vectorial closedness conditions.\\

{\bf Proposition  4.7.}  \ {\sl Let $E$ satisfy $\lambda \,E\subset E$ for all $\lambda >0$ and $\varphi_E^Q(y) \in {\bf R}$.

{\rm (a)} \  Let $Q-E$ be vectorially  closed and $\varphi_E^Q(y) >0$. Then $y\,\in\, \varphi_E^Q(y)\, Q-E$.

{\rm (b)} \  Let $Q+E$ be vectorially closed and $\varphi_E^Q(y) <0$. Then $y\,\in\,\varphi_E^Q(y)\,Q-E.$

{\rm (c)} \   Let $Q$ be convex and $\varphi_E^Q(y) =0$. Then $y\,\in\,-{\rm cl}_Q(E)$.}\\

{\sl  Proof.} \  (a) \  If $y\,\in\, \varphi_E^Q(y)\, Q-E$, then the result already holds. Now, assume that there exists a sequence $\{\epsilon_n\}$ with all $\epsilon_n >0$ and $\epsilon_n \rightarrow 0$ such that
$$y\,\in\,(\varphi_E^Q(y) +\epsilon_n)\, Q -E.$$
Thus,
$$\frac{y}{\varphi_E^Q(y) + \epsilon_n}\,\in\, Q-E.$$
Letting  $n\rightarrow \infty$, we have
$$\frac{y}{\varphi_E^Q(y)}\,\in\,{\rm vcl}(Q-E)\,=\, Q-E.$$
Thus, $y\,\in\,\varphi_E^Q(y)\, Q-E.$

(b) \  If $y\,\in\,\varphi_E^Q(y)\,Q-E$, then the result already holds. Now, assume that there exists a  sequence $\{\epsilon_n\}$ with all $\epsilon_n >0$ and $\epsilon_n \rightarrow 0$ such that
$$\varphi_E^Q(y) +\epsilon_n <0\ \ {\rm and} \ \  y\,\in\,(\varphi_E^Q(y) +\epsilon_n)\,Q-E.$$
Thus,
$$\frac{y}{\varphi_E^Q(y) + \epsilon_n}\,\in\,Q+E.$$
Letting $n\rightarrow\infty$, we  have
$$\frac{y}{\varphi_E^Q(y)}\,\in\, {\rm vcl}(Q+E)\,=\, Q+E.$$
Thus,
$y\,\in\,\varphi_E^Q(y)\, Q-E.$

(c) \  If $y\,\in\, \varphi_E^Q(y)\,Q-E\,=\,0-E\,=\, -E$, certainly, $y\,\in\,-{\rm cl}_Q(E)$.
Now, assume that there exists a sequence $\{\epsilon_n\}$ with all $\epsilon_n >0$ and $\epsilon_n \rightarrow 0$  such that
$$y\,\in\,\epsilon_n Q-E\,=\,-(E-\epsilon_n Q),\ \ \forall n.$$
Thus,
$$y\,\in\, \bigcap\limits_{n=1}^{\infty} -(E-\epsilon_n Q)\,\subset\,\bigcap\limits_{n=1}^{\infty}-(E-\epsilon_n Q_0)\,=\, -{\rm cl}_Q(E).$$\hfill\framebox[2mm]{}\\

{\bf Proposition   4.8.}  \ {\sl Let $Y$ be a topological vector space, $E$ be a closed convex cone and $Q$ be a compact convex set. Then $y\,\in\, \varphi_E^Q(y)\, Q - E$.}\\

{\sl  Proof.} \   Since $E$ is closed and $Q$ is compact, both $Q-E$ and $Q+E$ are closed, certainly are vectorially closed.
Now, applying Proposition 4.7 (a) and (b),  we have
$$y\,\in\,\varphi_E^Q(y)\,Q-E \ \ \   {\rm when} \ \  \  \varphi_E^Q(y) >0\ \  {\rm or} \ \  \varphi_E^Q(y) <0.$$
Next, we consider the case of $\varphi_E^Q(y) =0$. Since $E$ is closed, it is also locally closed and hence ${\rm cl}_Q(E)\,=\,E$.
Applying Proposition 4.7 (c), we have
$$y\,\in\,-{\rm cl}_Q(E)\,=\,-E\,=\,\varphi_E^Q(y) -E.$$\hfill\framebox[2mm]{}\\

Particularly, if $Y$ is a locally convex space and $E$ is a closed convex cone, the same result, i.e., $y\,\in\,\varphi_E^Q(y)-E$, still holds, even we only assume that $Q$ is weakly compact.

At the end of this section, we shall present a separation result on generalized nonconvex separation functional  in the setting of real vector spaces, which is similar to [18, Theorem 2.3.6]. For this, we need the notions of  algebraic interior, i.e., ${\rm cor}(E)$, and relatively algebraic interior, i.e., ${\rm icr}(E)$,  of set $E$, please refer to [8, 18, 27, 36, 53]. Here, we state the related notions in a slightly different way. Let $Y$ be a real vector space and $E\subset Y$ be nonempty. A point $y\in E$ is said to be a quasi-core point of $E$, denoted by $y\in {\rm qcor}(E)$, iff for any $v\in Y$ and  for any $\delta >0$, there exists  $0<\epsilon <\delta$ such that $y+\epsilon\, v\,\in\, E$. Moreover, let $B\subset Y$ be nonempty. A point $y\in E$ is said to be a quasi-core point of $E$ with respect to $B$, denoted by $y\in {\rm qcor}_B(E)$, iff for any $b\in B$ and  any $\delta >0$, there exists $0<\epsilon <\delta$ such that $y+\epsilon\, b\in E$.\\

{\bf Proposition  4.9.}  \ {\sl Let $E\subset Y$ satisfy $\lambda\, E\subset E$ for all $\lambda >0$ and let $Q\subset Y$ be nonempty. Suppose that the following conditions are satisfied:

{\rm (i)} \  $E+ (0, \infty)\cdot Q\, \subset\, {\rm qcor}_{-Q}(E).$

{\rm (ii)} \ $Q\cap {\rm qcor}(E)\, \not=\,\emptyset.$

{\rm (iii)} \  $ 0\not\in {\rm vcl}(Q+E).$

Then, for every  $y\in Y$, $\varphi^Q_E(y)\in {\bf R}$ and for any $\lambda\in {\bf R}$,
$$\{y\in Y:\, \varphi(y) <\lambda\}\,=\,\lambda\,Q- {\rm qcor}_{-Q}(E).$$
Particularly, if $\lambda =0$, then $\{y\in Y: \varphi(y) <0\}  = -{\rm qcor}_{-Q}(E)$. If a set $A\subset Y$ satisfies $A\cap -{\rm qcor}_{-Q}(E)\,=\, \emptyset$, then
$$\varphi^Q_E(-y) <0,\ \forall y\in {\rm qcor}_{-Q}(E) \ \ \  {\rm and} \ \  \  \varphi^Q_E(a) \geq 0,\ \forall a\in A.$$}\\

{\sl  Proof.} \  By condition (ii),  there exists $q\in Q\cap {\rm qcor}(E)$. For any $y\in Y$, there exists $\epsilon >0$ such that $q-\epsilon\,y \in E$.  From this,
$$y\,\in\, \frac{1}{\epsilon} q -E\,\subset\, \frac{1}{\epsilon} Q-E.$$ Thus,
$$\varphi^Q_E(y)\,\leq\,\frac{1}{\epsilon}\,<\,+\infty.$$
By condition (iii) and Proposition 4.2, we know that $\varphi^Q_E(y) \not= -\infty$.  Thus, we have
$\varphi^Q_E(y)\,\in\, {\bf R},\ \forall y\in Y.$ Next, we show that
$$\{y\in Y:\, \varphi^Q_E(y) < \lambda\}\,=\,\lambda\, Q- {\rm qcor}_{-Q}(E).$$
Assume that $\varphi^Q_E(y) < \lambda$. Then, there exists $\epsilon >0$ such that $\varphi^Q_E(y)  \leq \lambda-\epsilon < \lambda$ and
\begin{eqnarray*}
y\,&\in&\,(\lambda-\epsilon)\, Q -E\\
&\subset&\, \lambda\, Q-\epsilon\,Q-E\\
&\subset&\,\lambda\,Q- {\rm qcor}_{-Q}(E),
\end{eqnarray*}
where we have used condition (i) in the last step.

Conversely, assume that $y\,\in\, \lambda\,Q-{\rm qcor}_{-Q}(E)$. Then, there exists $q\in Q$ such that  $y\,\in\, \lambda\,q -{\rm qcor}_{-Q}(E).$
From this, $\lambda\,q-y\,\in\,{\rm qcor}_{-Q}(E)$. Hence, there exists $\epsilon >0$  such that
$\lambda\,q - y - \epsilon\,q\,\in\, E$. Thus,
$$y\,\in\,(\lambda-\epsilon) q -E\,\subset\,(\lambda-\epsilon)Q-E.$$
So, $\varphi^Q_E(y)\,\leq\,\lambda-\epsilon \,<\,\lambda.$\hfill\framebox[2mm]{}\\

\section*{\large\bf  5  \    EVPs with  set-valued objective functions and set-valued perturbations  }

\hspace*{\parindent}
In this section,  we assume that $(X, d)$ is a metric space, $Y$ is a real vector space quasi-ordered by a convex cone $K$, $H\subset K$ is a $K$-convex set and $f:\, X\rightarrow  2^Y$ is a set-valued map. We define a binary relation $\preceq$ on $X$ as follows: for any $x, x^{\prime}\in X$,
$$x^{\prime} \preceq x\ \ \   \Longleftrightarrow\ \  \  f(x) \subset  f(x^{\prime}) + d(x, x^{\prime})\, H +K.$$
Obviously, $x\preceq x$ for all $x\in X$. Thus, $\preceq$ satisfies the reflexive property.
Now, assume that $x^{\prime}\preceq  x$ and  $x^{\prime\prime}\preceq x^{\prime}$. Then,
$$f(x)\,\subset\, f(x^{\prime}) + d(x, x^{\prime})\, H +K \ \ \ \  {\rm  and}\ \ \ \
f(x^{\prime})\,\subset\, f(x^{\prime\prime}) + d(x^{\prime}, x^{\prime\prime})\, H +K.$$
From the above two inclusions, we have
\begin{eqnarray*}
f(x) \,&\subset&\, f(x^{\prime\prime})  + d(x^{\prime}, x^{\prime\prime})\,H+ d(x, x^{\prime})\, H +K\\
&\subset&\,f(x^{\prime\prime})  + d(x^{\prime}, x^{\prime\prime})(H+K)+ d(x, x^{\prime})(H+K) +K\\
&=&\, f(x^{\prime\prime}) +(d(x, x^{\prime}) + d(x^{\prime}, x^{\prime\prime}))\,(H+K) +K\\\
&\subset&\, f(x^{\prime\prime}) + d(x, x^{\prime\prime})\, H +K.
\end{eqnarray*}
That is,  $x^{\prime\prime}\preceq x$. Thus, $\preceq$ satisfies the transitive property. Hence $\preceq$ is a quasi-order.
For any $x\in X$, put
$$S(x) :=\,\{x^{\prime}\in X:\, x^{\prime}\preceq x\}.$$
As $\preceq$ is a quasi-order, we have
$x\in S(x)$ for all $x\in X$; and $S(x^{\prime})\,\subset\, S(x)$ for all $x^{\prime}\in S(x)$.\\

{\bf Definition  5.1} (See [42]). \  Let $S(\cdot):\, X\rightarrow  2^X\backslash\{\emptyset\}$ be a set-valued map satisfying  $S(x^{\prime})\subset S(x)$ for all $x^{\prime}\in S(x)$.
$S(\cdot)$ is said to be dynamically closed at $x\in X$ if $(x_n)\subset  S(x)$ such that $x_{n+1}\in S(x_n)$ for all  $n$ and $x_n\rightarrow \bar{x}$ then $\bar{x}\in S(x)$. In this case, we also say that $S(x)$ is dynamically closed. Moreover, $(X,d)$ is said to be $S(x)$-dynamically complete if for any Cauchy sequence $(x_n)\subset  S(x)$ such that $x_{n+1}\in S(x_n)$ for all $n$, there exists $\bar{x}\in X$ such that $x_n\rightarrow \bar{x}$.\\

We remark that a property similar  to the above dynamical closedness, i.e., so-called the limit monotonicity property, was also introduced by  Bao and Mordukhovich (See [4, 5]).
Moreover, some useful  notions, for example,  $\preceq$-completeness and $\preceq$-lower closedness on quasi-order $\preceq$,  were introduced  by Khanh and Quy (see [29]).
A quasi-order $\preceq$ is said to be lower closed if for any sequence $(x_n)\subset X$  such that $x_{n+1}\preceq x_n$ for all $n$ and converging to $\bar{x}$, one has $\bar{x}\preceq x_n$ for all $n$.
$S(x)$ is said to be $\preceq$-complete if every Cauchy sequence $(x_n)\subset S(x)$ such that $x_{n+1}\in S(x_n)$ for all $n$, is convergent to a point of $S(x)$.
Obviously, $\preceq$ being lower closed implies that $S(x)$ is dynamically closed for all $x\in X$. And $S(x)$ being $\preceq$-complete implies that $(X, d)$ is $S(x)$-dynamically complete.
By using lower closed quasi-order, under very weak conditions Khanh and Quy [29] obtained the following very general set-valued EVP, where the perturbation consists of a convex subset of the ordering cone multiplied by the distance function.\\

{\bf Theorem  5.A} (See [29, Theorem 3.2]).  \ {\sl Let $(X, d)$ be a metric space, $Y$ be a locally convex space, $K$ be a convex cone in $Y$ and $H\subset  K$ be a convex set such that $0\not\in {\rm cl}(H+K),$
$x_0\in X$ and $f:\, X\rightarrow 2^Y$ be a set-valued map.
Suppose that the following conditions hold

{\rm (i)} \  $S(x_0)$ is $\preceq$-complete;

{\rm (ii)} \ $f(S(x_0))$ is quasi $K$-lower bounded (i.e., quasibounded from below);

{\rm (iii)} \  $\preceq$ is lower closed.

Then, there exists $\bar{x}\in X$ such that

{\rm (a)} \  $f(x_0)\,\subset\, f(\bar{x}) + d(x_0, \bar{x})\,H +K;$

{\rm (b)} \   $f(\bar{x})\,\not\subset\,  f(x) + d( \bar{x}, x)\, H +K, \ \ \ \forall x\in X\backslash\{\bar{x}\}.$
}\\

In order to deduce our set-valued EVPs, we need the following lemma, which is indeed a corollary of the pre-order principle in [44].\\

{\bf Lemma  5.2} (See [44, Theorem 3.1]). \ {\sl Let $X$ be a nonempty set,  $Y$ be a real vector space, $K\subset Y$  be a convex cone specifying a quasi-order  $\leq_K$ on  $Y$, $f:\,  X\rightarrow  2^Y$ be a set-valued map
and $F_{\lambda}:\, X\times X \rightarrow 2^K\backslash\{\emptyset\},\ \lambda\in\Lambda,$ be a family of set-valued bimaps satisfying the triangle inequality property {\rm (see [16])}, i.e., for each $x_i\in X,\ i=1,2,3,$ and $\lambda\in\Lambda$ there exists $\mu,\nu\in\Lambda$ such that
$$F_{\mu}(x_1, x_2)+ F_{\nu}(x_2, x_3) \,\subset\, F_{\lambda}(x_1, x_3) +K.$$
Let $x_0\in X$ such that
$$S(x_0):=\,\{x\in X:\, f(x_0)\subset f(x) + F_{\lambda}(x, x_0) + K,\ \forall \lambda\in\Lambda\}\,\not=\,\emptyset.$$
Suppose that there exists a $K$-monotone extended real function $\xi:\, {\bf R}\cup\{\pm\infty\}$ satisfying the following assumptions:

{\rm (D)} \  $-\infty < \inf\,\xi\circ f(S(x_0)) < +\infty;$

{\rm (E)} \ for any $x\in S(x_0)$ with $-\infty < \xi\circ f(x) <+\infty$ and for any $x^{\prime}\in S(x)\backslash \{x\}$, one has $\inf\,\xi\circ f(x) > \inf\,\xi\circ  f(x^{\prime});$

{\rm (F)} \  for any sequence $(x_n)\subset S(x_0)$ with $x_n\in S(x_{n-1}),\ \forall n,$ such that $\inf\,\xi\circ f(x_n) -\inf\,\xi\circ  f(S(x_{n-1}))\,\rightarrow\, 0,$ there exists $u\in X$ such that $u\in S(x_n),\ \forall n.$

Then, there exists $\bar{x}\in X$ such that

{\rm (a)} \ $f(x_0)\,\subset\, f(\bar{x}) + F_{\lambda}(\bar{x}, x_0) +K,\ \forall \lambda\in\Lambda;$

{\rm (b)} \  $\forall x\in X\backslash\{\bar{x}\},\ \exists \lambda\in\Lambda$ such that $f(\bar{x})\,\not\subset\, f(x) + F_{\lambda}(x, \bar{x}) + K.$
}\\

By checking the proofs of [44, Theorem 2.1] and [44, Theorem 3.1], we see that ``for any $x\in S(x_0) \cdots$" in assumption (E) of Lemma 5.2 can be replaced by ``for any $x\in S(x_0)\backslash\{x_0\}\cdots$".
Using Lemma 5.2 and generalized nonconvex separation functionals we shall give a general version of  EVP with set-valued objective function and set-valued perturbation,
 which improves Theorem 5.A by weakening some conditions (See Remark 5.4).\\

{\bf Theorem 5.3. }  \ {\sl Let $(X, d)$ be a metric space, $Y$ be a real vector space, $K\subset Y$ be a convex cone, $H\subset K$ be a $K$-convex set with $0\not\in {\rm vcl}(H+K)$,  $f:\, X\rightarrow 2^Y$ be a set-valued map and $x_0\in X$. For any $x\in X$, put
$$S(x):=\,\{x^{\prime}\in X:\, f(x)\,\subset\, f(x^{\prime})+ d(x, x^{\prime})H +K\}.$$
Suppose that the following conditions are satisfied:

{\rm (i)} \  $(X,d)$ is $S(x_0)$-dynamically complete;

{\rm (ii)} \  there exists  $\epsilon >0$ such that $f(x_0)\,\not\subset\, f(S(x_0)) + \epsilon H +K;$

{\rm (iii)} \  for any $x\in S(x_0)$, $S(x)$ is dynamically closed.

Then, there exists $\bar{x}\in X$ such that

{\rm (a)} \  $f(x_0)\,\subset\, f(\bar{x}) + d(x_0, \bar{x}) H +K;$

{\rm (b)} \  $f(\bar{x})\,\not\subset\, f(x) + d(\bar{x}, x) H +K,\ \forall x\in X\backslash\{\bar{x}\}.$
}\\

{\sl Proof.} \ Define $F:\, X\times X\rightarrow 2^K\backslash\{\emptyset\}$ as follows
$$F(x, x^{\prime}): =\, d(x, x^{\prime})\, H,\ \  \forall x,\, x^{\prime}\in X.$$
Since $H$ is a $K$-convex set, for any $x_1, x_2, x_3\in X$, we have
\begin{eqnarray*}
&\ &\,d(x_1, x_2)\, H + d(x_2, x_3) \,H\\
&\subset&\,d(x_1, x_2)\, (H+K) + d(x_2, x_3) \,(H+K)\\
&=&\, (d(x_1, x_2) + d(x_2 , x_3))\,(H+K)\\
&=&\,(d(x_1, x_2) + d(x_2, x_3)) H + K\\
&\subset&\, d(x_1, x_3)\, H +K.
\end{eqnarray*}
Thus, a singleton $\{F\}$ satisfies the triangle inequality property.
We shall apply Lemma 5.2 to prove the result.

By condition (ii), there exists $y_0\in f(x_0)$ such that
$$y_0\,\not\in\, f(S(x_0)) + \epsilon H +K.\eqno{(5.1)}$$
Define $\xi:\, Y\rightarrow {\bf R}\cup\{\pm\infty\}$ as follows:
$$\xi(y):=\,\varphi^H_K(y-y_0),\ \ \  y\in Y.$$
By Proposition 4.6(b),  we know that $\xi$ is a $K$-monotone extended real function. We shall show that assumptions (D), (E), (F) in Lemma 5.2 are  all satisfied.

{\bf Step 1}  \  Show that (D) is satisfied. \  By (5.1)
$$(f(S(x_0))-y_0)\cap (-\epsilon H -K)\,=\,\emptyset,$$
hence
$$\xi(y) = \varphi^H_K(y-y_0)\,\geq\, -\epsilon, \ \ \ \forall y\in f(S(x_0)).$$
Also, $x_0\in S(x_0)$, $y_0\in f(x_0)\subset f(S(x_0))$ and $\xi(y_0) = \varphi^H_K(y_0-y_0) =\varphi^H_K(0) =0$.
So $\inf\,\xi\circ f(S(x_0))\,\leq\, \inf\, \xi\circ f(x_0) \leq \xi(y_0) = 0.$
Therefore,
$$-\infty < -\epsilon \leq \inf\,\xi\circ f(S(x_0)) \leq 0 < +\infty.$$
Thus, (D) is satisfied.

{\bf Step 2}  \   Show that (E) is satisfied. \  Take any $x\in S(x_0)\backslash\{x_0\}$ with $-\infty <\inf\,\xi\circ f(x) < +\infty$. Since $f(x_0)\,\subset\, f(x) + d(x_0, x) H +K$, we have $y_0\,\in\, f(x) + d(x_0, x)\, H +K.$
For all the elements in $f(x)$, we divide them into the following two sets:  the set $\{y\in f(x):\, y_0\,\in\, y+d(x_0, x)\, H +K\}$ denoted by $f(x)_1$; the set $\{y\in f(x):\, y_0\,\not\in\, y+ d(x_0, x)\, H +K\}$ denoted by $f(x)_2$.
Obviously,
$$f(x) =  f(x)_1 \cup f(x)_2,\ \  f(x)_1\cap f(x)_2 = \emptyset, \ \  {\rm and} \ \   f(x)_1\not= \emptyset.$$
For $y\in  f(x)_1$, we have $$y_0\,\in\, y+ d(x_0, x)\,H +K\ \ \  {\rm and}\ \ \  y-y_0\,\in\, -d(x_0 , x)\, H -K.$$ Thus,
$$\xi(y):= \varphi^H_K(y-y_0) \leq  -d(x_0, x) < 0, \ \ \  \forall y\in f(x)_1.\eqno{(5.2)}$$
For $y\in f(x)_2$, we have
$$y_0 \not\in y +d(x_0 , x)\,H +K  \ \ \ {\rm and} \ \ \   y-y_0 \not\in  -d(x_0, x)\, H -K.$$ Thus,
$$\xi(y):= \varphi^H_K(y-y_0) \geq  -d(x_0, x), \ \ \  \forall y\in f(x)_2.\eqno{(5.3)}$$
Combining (5.2) and (5.3) and remarking $f(x)_1\not=\emptyset$, we have
$$\inf\,\xi\circ f(x)  = \inf \xi\circ f(x)_1 \leq -d(x_0, x) < 0.\eqno{(5.4)}$$
Take any $x^{\prime} \in S(x)\backslash \{x\}$. Then $x\not= x^{\prime}$ and
$$ f(x)\,\subset\, f(x^{\prime}) + d(x, x^{\prime})\, H +K.$$
For any $y\in f(x)_1$, there exists $y^{\prime}\in f(x^{\prime})$ such that
$$y \in y^{\prime} +d(x, x^{\prime})\, H +K\ \ \ {\rm and}\ \ \  y^{\prime}-y \in -d(x, x^{\prime})\, H -K.$$
Thus,
$$\varphi^H_K(y^{\prime}-y)\,\leq\, -d(x, x^{\prime}) \,<\, 0.\eqno{(5.5)}$$
By (5.2) and (5.5), and using Proposition 4.6(d), we have
$$\varphi^H_K(y^{\prime}-y_0)\,\leq\, \varphi^H_K(y^{\prime} -y)  + \varphi^H_K(y-y_0). \eqno{(5.6)}$$
From (5.5) and (5.6), we have
\begin{eqnarray*}
d(x, x^{\prime})\,&\leq&\, -\varphi^H_K(y^{\prime}-y)\\
&\leq&\,\varphi^H_K(y-y_0) - \varphi^H_K(y^{\prime}-y_0)\\
&\leq&\, \varphi^H_K(y-y_0) -\inf_{y^{\prime}\in f(x^{\prime})}\varphi^H_K(y^{\prime} -y_0)\\
&=&\,\varphi^H_K(y-y_0) -\inf\,\xi\circ f(x^{\prime}).
\end{eqnarray*}
The above inequality holds for all $y\in f(x)_1$.
Thus,
\begin{eqnarray*}
d(x, x^{\prime}) \,&\leq&\, \inf_{y\in f(x)_1} \varphi^H_K(y-y_0) -\inf\,\xi\circ f(x^{\prime})\\
&=&\, \inf\,\xi\circ f(x)_1 - \inf\,\xi\circ f(x^{\prime})\\
&=&\, \inf\,\xi\circ f(x) - \inf\, \xi\circ  f(x^{\prime}),
\end{eqnarray*}
where the last equality is due to (5.4).
From this, we have$$\inf\,\xi\circ f(x)\,\geq\, \inf\,\xi\circ f(x^{\prime}) + d(x, x^{\prime}) \,>\, \inf\,\xi\circ f(x^{\prime}).$$
That is, (E) is satisfied.

{\bf Step 3}  \  Show that (F) is satisfied. \  Let a sequence $(x_n)\subset  S(x_0)$ with $x_n\in S(x_{n-1}),\ \forall n,$ satisfy
$$\inf\,\xi\circ f(x_n) - \inf\,\xi\circ f(S(x_{n-1})) \,\rightarrow\, 0 \ \ (n\rightarrow\infty).$$
We may take a positive sequence $(\epsilon_n)$ convergent to $0$ such that for every $n$,
$$\inf\,\xi\circ f(x_n) - \inf\,\xi\circ f(S(x_{n-1}))\,<\,\epsilon_n.\eqno{(5.7)}$$
We shall show that there exists $u\in X$ such that $u\in S(x_n),\ \forall n.$

If there exists a sequence $n_1 < n_2 < n_3 < \cdots$ such that $x_{n_1} = x_{n_2} = \cdots = x_{n_i} = \cdots$, then we may take $u$ to be the common element $x_{n_i}$.  Obviously, $u\in S(x_n),\ \forall n$.
Hence, without loss of generality, we may assume that $x_n \not= x_m$ if $n\not= m$. For each $n\in {\bf N}$, $x_n\in S(x_0)$, so
$$y_0\,\in\, f(x_0)\,\subset\, f(x_n) + d(x_0, x_n)\, H +K.$$
Hence, there exists $y_n^{\prime}\in f(x_n)$ such that
$$y_0\,\in\, y_n^{\prime} + d(x_0, x_n) \, H +K, \ \ \  {\rm i.e.,} \ \ \  y_n^{\prime} - y_0\,\in\, -d(x_0, x_n)\, H - K.$$
Thus,
$$\varphi^H_K(y_n^{\prime}-y_0) \leq -d(x_0, x_n) < -\frac{1}{2} d(x_0, x_n) < 0.$$
So $$\inf\,\{\varphi^H_K(y_n^{\prime\prime} - y_0):\, y_n^{\prime\prime} \in f(x_n)\}\,\leq\, \varphi^H_K(y_n^{\prime}-y_0)\,<\, -\frac{1}{2} d(x_0, x_n)  \,<\,0.$$
That is, $$\inf\, \xi\circ f(x_n) \,<\, -\frac{1}{2} d(x_0, x_n) \,<\, 0.$$
Combining this with (5.7), for each $n$, we may take
 $y_n\in f(x_n)$  such that
$$\xi(y_n) \, <\, \min\{\inf\,\xi\circ f(S(x_{n-1})) +\epsilon_n,\, -\frac{1}{2} d(x_0, x_n)\}. \eqno{(5.8)}$$
For $m>n$, $x_m\in S(x_n)$. So
$$y_n\,\in\, f(x_n)\,\subset\, f(x_m) + d(x_n, x_m) \,H +K.$$
Thus, there exists $y_{mn}\in f(x_m)$ such that
$$y_n\,\in\,y_{mn} + d(x_n, x_m)\,H +K,\ \ \ {\rm i.e.,}\ \ \ y_{mn}-y_n\,\in\, -d(x_n, x_m)\, H-K.$$
Hence,
$$\varphi^H_K(y_{mn}-y_n)\,\leq\,-d(x_n, x_m)\eqno{(5.9)}$$
and
$$d(x_n, x_m)\,\leq\,-\varphi^H_K(y_{mn}-y_n).\eqno{(5.10)}$$
By (5.8) and (5.9), we know that
$$\varphi^H_K(y_n-y_0) = \xi(y_n) < -\frac{1}{2} d(x_0, x_n) < 0\ \ \ {\rm  and} \ \ \ \varphi^H_K(y_{mn}-y_n) < 0.$$
Thus, we can apply Proposition 4.6(d) and have
$$\varphi^H_K(y_{mn}- y_0)\,\leq\, \varphi^H_K(y_{mn}-y_n) + \varphi^H_K(y_n -y_0).$$
From this,
$$-\varphi^H_K(y_{mn} -y_n)\,\leq\,\varphi^H_K(y_n-y_0) -\varphi^H_K(y_{mn} -y_0). \eqno{(5.11)}$$
Combining (5.10), (5.11) and (5.8), we have
\begin{eqnarray*}
d(x_n, x_m)\,&\leq&\,-\varphi^H_K(y_{mn}-y_n)\\
&\leq&\,\varphi^H_K(y_n - y_0) - \varphi^H_K(y_{mn} -y_0)\\
&=&\, \xi(y_n) - \xi(y_{mn})\\
&<&\, \inf\,\xi\circ f(S(x_{n-1})) +\epsilon_n  -\inf\,\xi\circ f(S(x_{n-1}))\\
&=&\,\epsilon_n,
\end{eqnarray*}
where we have used  $y_{mn} \in f(x_m) \subset f(S(x_{n-1}))$.
Since $\epsilon_n \rightarrow 0$, we know that $(x_n)\subset S(x_0)$ is a Cauchy sequence. Remark that $x_n\in S(x_{n-1})$. By condition (i), there exists  $u\in X$ such that $x_n \rightarrow u\ (n\rightarrow\infty)$. On the other hand, $(x_{n+p})_{p\in {\bf N}}\subset S(x_n),\ x_{n+p} \in S(x_{n+p-1})$ and $x_{n+p} \rightarrow u\ (p\rightarrow\infty)$. By condition (iii), we conclude that $u\in S(x_n)$. Thus,  (F) is satisfied.

Now, we can apply Lemma 5.2 and obtain the result.\hfill\framebox[2mm]{}\\

{\bf Remark 5.4.} \  Comparing Theorem 5.A with Theorem 5.3, we see that $Y$ being a locally convex space is replaced by $Y$ being a real vector space; $0\not\in {\rm cl}(H+K)$ is replaced by $0\not\in {\rm vcl}(H+K)$; and $H\subset K$ being a convex set is replaced by $H\subset K$ being a $K$-convex set. Besides, condition (ii) in Theorem 5.A ``$f(S(x_0))$ is quasi $K$-lower bounded" is replaced by condition (ii) in Theorem 5.3 ``There exists  $\epsilon >0$ such that $f(x_0)\not\subset f(S(x_0)) +\epsilon\, H+K$". Obviously, the latter is strictly weaker than the former (For details, see Section 3).\\

If condition (iii) in Theorem 5.3 is replaced by a stronger condition: for any $x\in S(x_0)$, $S(x)$ is closed, then we can obtain the following corollary, which improves [29, Corollary 3.4]  and [51, Theorem 4.2].\\

{\bf Corollary   5.5.}  \  {\sl Let $(X, d)$ be a metric space, $Y$ be a real vector space, $K\subset Y$ be a convex cone, $H\subset K$ be a $K$-convex set with $0\not\in {\rm vcl}(H+K)$,  $f:\, X\rightarrow 2^Y$ be a set-valued map, $x_0\in X$, $\epsilon >0$ and $\lambda >0$. For any $x\in X$, put
$$S(x):=\,\{x^{\prime}\in X:\, f(x)\,\subset\, f(x^{\prime})+ \frac{\epsilon}{\lambda} d(x, x^{\prime})H +K\}.$$
Suppose that the following conditions are satisfied:

{\rm (i)} \  $(X,d)$ is $S(x_0)$-dynamically complete;

{\rm (ii)} \  $f(x_0)\,\not\subset\, f(S(x_0)) + \epsilon H +K;$

{\rm (iii)} \  for any $x\in S(x_0)$, $S(x)$ is  closed (or, $S(x)$ is dynamically closed).

Then, there exists $\bar{x}\in X$ such that

{\rm (a)} \  $f(x_0)\,\subset\, f(\bar{x}) + ({\epsilon}/{\lambda}) d(x_0, \bar{x}) H +K;$

{\rm (b)} \  $f(\bar{x})\,\not\subset\, f(x) + ({\epsilon}/{\lambda}) d(\bar{x}, x) H +K,\ \forall x\in X\backslash\{\bar{x}\};$

{\rm (c)} \  $d(x_0, \bar{x}) \leq \lambda$.
}\\

{\sl Proof.} \   Obviously, $S(x)$ being closed implies that $S(x)$ is dynamically closed. Put
$$\tilde{d}(x, x^{\prime}) \,=\, \frac{\epsilon}{\lambda}\, d(x, x^{\prime}),\ \ \ \forall x, x^{\prime} \in X.$$
Replacing $d$ by $\tilde{d}$ in Theorem 5.3, we conclude that there exists $\bar{x}\in X$  such that  (a) and (b)  hold. If $d(x_0, \bar{x})  > \lambda$, then by (a),
\begin{eqnarray*}
f(x_0)\,&\subset&\, f(\bar{x}) + \frac{\epsilon}{\lambda}\, d(x_0, \bar{x})\,H +K\\
&\subset&\, f(\bar{x}) + \frac{\epsilon}{\lambda}\, {\lambda}\,H +K\\
&=&\, f(\bar{x}) + \epsilon\, H +K,
\end{eqnarray*}
where $\bar{x}\in S(x_0)$.
This contradicts condition (ii). Thus, (c) holds.\hfill\framebox[2mm]{}\\

A set-valued map $f:\, X\rightarrow 2^Y$ is said to have $K_H$-closed values iff for any $x\in X$ and any $\alpha >0$, $f(x) +\alpha\, H +K$ is closed (see [29, 44]).
As we have seen at the beginning of  Section 4, vectorial closedness and $v_0$-closedness are weaker than topological closedness.
Let $v_0\in Y$. Then
$f$ is said to have $K_H$-$v_0$-closed values iff for any $x\in X$ and any $\alpha >0$, $f(x) + \alpha\,H +K$ is $v_0$-closed. Remark that every closed set in $Y$ is $v_0$-closed for all $v_0\in Y$, but the converse is not true.
So, $f$ having $K_H$-closed values implies that $f$ has $K_H$-$v_0$-values for all $v_0\in Y$, but the converse is not true.

The following corollary extends and improves [29, Corollary 3.5] and [44, Theorem 4.2$^{\prime}$].\\

{\bf  Corollary   5.6.} \ {\sl Let $(X, d),\, Y,\, K,\,H,\, f,\, x_0,\, \epsilon >0,\,\lambda >0, $ and $S(x)$ be the same as in Corollary 5.5. Suppose that the following conditions are satisfied:

{\rm (i)} \  $(X,d)$ is $S(x_0)$-dynamically complete;

{\rm (ii)} \   $f(x_0)\,\not\subset\, f(S(x_0)) + \epsilon H +K;$

{\rm (iii)} \  $f$ is $K$-slm and has $K_H$-$v_0$-closed values, where $v_0\in K+H$.

Then, there exists $\bar{x}\in X$ such that {\rm (a)}, {\rm (b)} and {\rm (c)} in Corollary 5.5 hold.}\\

{\sl Proof.} \   By Theorem 5.3 we only need to prove that for any $x\in S(x_0)$, $S(x)$ is dynamically closed. Let $(x_n)\subset S(x)$ satisfy that  $x_{n+1}\in  S(x_n),\ \forall n$, and $x_n\rightarrow \bar{x}$.
Let $n\in {\bf N}$ be given. For every  $k\in {\bf N}$, $x_{n+k} \in S(x_n)$,  that is,
$$f(x_n)\,\subset\, f(x_{n+k}) +\frac{\epsilon}{\lambda} d(x_n, x_{n+k})\, H +K.\eqno{(5.12)}$$
Since $f(x_n) \,\subset\, f(x_{n+1}) +K,\ \forall n$, and $x_n\rightarrow\bar{x}$, by the condition that $f$ is $K$-slm, we have  $f(x_n)\,\subset\, f(\bar{x}) +K.$ Combining this with (5.12) we have
$$f(x_n)\,\subset\, f(\bar{x})  +\frac{\epsilon}{\lambda} d(x_n, x_{n+k})\, H + K.\eqno{(5.13)}$$
Next, we consider the following two different cases.

{\bf Case 1.} \  There exists $k\in {\bf N}$ such that $d(x_n, x_{n+k}) \geq\, d(x_n, \bar{x}).$
Combining this with (5.13) we have
$$f(x_n)\,\subset\, f(\bar{x}) + \frac{\epsilon}{\lambda} d(x_n, \bar{x})\, H +K,\ \ {\rm i.e.,} \ \  \bar{x} \in S(x_n).$$

{\bf Case 2.} \  For all $k\in {\bf N}$, $d(x_n, x_{n+k}) \,<\, d(x_n, \bar{x}).$
Thus, from  (5.13) we have
\begin{eqnarray*}
&\  &\, f(x_n)  + \frac{\epsilon}{\lambda} (d(x_n, \bar{x}) - d(x_n, x_{n+k})) v_0\\
&\subset&\, f(\bar{x}) +\frac{\epsilon}{\lambda} d(x_n, x_{n+k}) H +K+ \frac{\epsilon}{\lambda} (d(x_n, \bar{x}) - d(x_n, x_{n+k})) v_0\\
&\subset&\,   f(\bar{x}) +\frac{\epsilon}{\lambda} d(x_n, x_{n+k})(H +K) + \frac{\epsilon}{\lambda} (d(x_n, \bar{x}) - d(x_n, x_{n+k}))(H+K) +K\\
&=&\,  f(\bar{x}) + \frac{\epsilon}{\lambda} d(x_n, \bar{x})  (H+K) +K\\
&=&\, f(\bar{x}) + \frac{\epsilon}{\lambda} d(x_n, \bar{x}) \,H +K.  \hspace*{8cm} (5.14)
\end{eqnarray*}
Since $d(x_n, \bar{x}) - d(x_n, x_{n+k})\,\rightarrow\, 0 \ (k\rightarrow\infty)$ and $f(\bar{x}) + (\epsilon/\lambda) d(x_n, \bar{x})\, H +K$ is $v_0$-closed, from (5.14) we have
$$f(x_n)\,\subset\, f(\bar{x}) + \frac{\epsilon}{\lambda} d(x_n, \bar{x})\, H +K,\ \ \  {\rm i.e.,} \ \  \   \bar{x} \in S(x_n).$$
Synthesizing Case 1 and Case 2, we conclude that $\bar{x} \in S(x_n) \subset S(x)$. Thus, we have shown that $S(x)$ is dynamically closed. \hfill\framebox[2mm]{}\\

From Corollary 5.6, we can also  deduce the following corollary, which extends and improves [24, Theorem 3.1] and [41, Theorem 3.1].\\

{\bf Corollary  5.7.} \  {\sl Let $(X, d),\, Y,\, K,\, H,\, f,\, x_0,\, \epsilon >0$ and  $\lambda >0$  be the same as in Corollary 5.5.  Suppose that the following conditions are satisfied:

{\rm (i)} \ $(X,d)$ is $(f, K)$-lower complete;

{\rm (ii)} \  $f(x_0) \not\subset f(X) + \epsilon\,H +K;$

{\rm (iii)} \ $f$ is $K$-slm and has $K_H$-$v_0$-closed values, where $v_0$ is a certain point in $H+K$.

Then, there exists $\bar{x}\in X$ such that {\rm  (a), (b)} and {\rm (c)} in Corollary 5.5 hold.}\\

{\sl Proof.} \  It is easy to see that $(X, d)$ being $(f, K)$-lower complete implies that $(X, d)$ is $S(x_0)$-dynamically complete. Now, applying Corollary 5.6 we immediately obtain the result.\hfill\framebox[2mm]{}

In particular, if $H$ is exactly a singleton $\{k_0\}\subset K\backslash -{\rm vcl}(K)$, then $f$ having $K_H$-closed values becomes an easier form, i.e.,  $f$ has $K$-closed values; and $f$ having $K_H$-$v_0$-closed values becomes that
$f$  has $K$-$v_0$-closed values, i.e., $f(x)+K$ is $v_0$-closed for all $x\in X$. Thus, by using  Corollary 5.6, we can obtain the following corollary, which improves [24, Theorem 3.1], [41, Theorem 3.1], [34, Theorem 3.1] and [44, Corollary 3.6].\\

{\bf Corollary  5.8.} \  {\sl Let  $(X, d)$ be a metric space, $Y$ be a real vector space,  $K\subset Y$ be a convex cone and $k_0\in K\backslash -{\rm vcl}(K)$. Let $f:\, X\rightarrow 2^Y$ be $K$-slm and have
$K$-$v_0$-closed values, where $v_0\in k_0 +K$.
Suppose that $\epsilon >0$ and $x_0\in X$ such that
$$f(x_0) \not\subset f(X) +\epsilon k_0 +K,$$ and suppose that $(X, d)$ is $S(x_0)$-dynamically complete (or, $(X, d)$ is $(f, K)$-lower complete), where $S(x_0) = \{x\in X:\,  \,f(x_0)\subset f(x) + (\epsilon/\lambda)d(x_0, x) k_0 +K\}$ and  $\lambda >0$ is a constant. Then, there exists  $\bar{x}\in X$ such that

{\rm (a)} \ $f(x_0)\subset f(\bar{x}) + (\epsilon/\lambda) \,d(x_0, \bar{x})\,k_0 + K;$

{\rm (b)} \ $\bar{x}$ is a strict minimizer of the map $x\mapsto f(x) + (\epsilon/\lambda)\, d(\bar{x}, x)\,k_0$, i.e.,
$$f(\bar{x})\,\not\subset\, f(x) + (\epsilon/\lambda)\, d(\bar{x}, x)\, k_0 +K,\ \ \forall x\in X\backslash\{\bar{x}\};$$

{\rm (c)} \  $d(x_0, \bar{x}) \leq \lambda.$}\\

\section*{\large\bf  6  \    EVPs for $\epsilon$-efficient  solutions in set-valued optimization  }

\hspace*{\parindent}
In this section, we always assume that $(X, d)$ is a metric space, $Y$ is a real vector space and $K\subset Y$ is a pointed convex cone specifying a partial order $\leq_K$ on $Y$.
Let us consider the following
vector optimization problem:
$$   {\rm Min}\{f(x):\, x\in S\}, \eqno{(6.1)}$$
where $f: X\rightarrow Y$ is a vector-valued map and $S$ is a
nonempty closed subset of $X$. A point $x_0\in S$ is called an
efficient solution of (6.1) if
$$  (f(S)-f(x_0))\cap (-K\backslash\{0\})=\emptyset,$$
where $f(S)$ denotes the set $\cup_{x\in S}\{f(x)\}$.

Guti\'{e}rrez, Jim\'{e}nez and Novo [20, 21] introduced the $(C, \epsilon)$-efficiency concept, which extends and unifies several $\epsilon$-efficiency notions.\\

{\bf Definition   6.1} ([21, 22]).  \   A nonempty set $C\subset Y$ is
coradiant if $\cup_{\beta \geq 1}\beta C = C$.\\

{\bf Definition 6.2} ([21, 22]).  \  Let $K$ be an ordering cone,
$C\subset K\backslash\{0\}$ be a coradiant set and let $\epsilon >0$. A
point $x_0\in S$ is a $(C,\epsilon)$-efficient solution of problem
{\rm (6.1)} if $(f(S)-f(x_0))\cap (-\epsilon C) =\emptyset.$
In this case, we also denote  $x_0\in AE(C, \epsilon)$.\\

In particular, if $C:= H+K$, where $H\subset K\backslash\{0\}$, then we can easily verify that $C$ is a coradiant set and $C\subset K\backslash\{0\}$. Thus, we obtain the concept of approximate efficiency due to N\'{e}meth.\\

{\bf Definition  6.3} ([22, 38]).  \  Let $H\subset
K\backslash\{0\}$ and let $\epsilon >0$. A point $x_0\in S$ is said to
be an $\epsilon$-efficient solution of {\rm (6.1)} in the sense of
N\'{e}meth (with respect to $H$) if $(f(S)-f(x_0))\cap (-\epsilon H-K)\,=\,\emptyset$.
In this case, we also denote $x_0\in AE(C_H, \epsilon)$, where $C_H = H+K$.\\

Next, let us consider the following set-valued optimization problem:
$$   {\rm Min}\{f(x):\, x\in S\}, \eqno{(6.2)}$$
where $f: X\rightarrow 2^Y$ is a set-valued map and $S$ is a
nonempty closed subset of $X$. A point $x_0\in S$ is called an
efficient solution of (6.2) if there exists $y_0\in f(x_0)$ such that
$$  (f(S)-y_0)\cap (-K\backslash\{0\})=\emptyset,$$
where $f(S)$ denotes the set $\cup_{x\in S}f(x)$.

Moreover, we  have the
 $(C, \epsilon)$-efficiency concept  in set-valued optimization.\\

{\bf Definition   6.4.} \  Let $C\subset K\backslash\{0\}$ be a coradiant set and let $\epsilon >0$.  A point $x_0\in S$ is called a $(C, \epsilon)$-efficient solution of problem (6.2)
if there exists $y_0\in f(x_0)$ such that
$$(f(S) - y_0)\cap (-\epsilon\, C)\,=\,\emptyset.$$
In this case, we also denote $x_0\in AE(C,\epsilon)$.\\

Similarly, we also have the  concept of approximate efficiency
due to N\'{e}meth in set-valued optimization, which extends the corresponding concept in vector optimization in [21, 22].\\

{\bf Definition  6.5.} \  Let $H\subset K\backslash\{0\}$ and $\epsilon >0$. A point $x_0\in S$ is said to be an $\epsilon$-efficient solution of (6.2) in the sense of N\'{e}meth (with respect to $H$) if there exists $y_0\in f(x_0)$ such that
$$(f(S) - y_0)\cap (-\epsilon\,H -K)\,=\, \emptyset.$$
In this case, we also denote $x_0\,\in\, AE(C_H, \epsilon)$, where $C_H \,=\, H+K.$\\

From Corollary 5.5 we can obtain a set-valued  EVP for approximately efficient solutions in set-valued optimization as follows.\\

{\bf Theorem 6.6. } \ {\sl Let $(X,d)$ be a metric space, $Y$ be a real vector space, $K\subset Y$ be a pointed convex cone specifying a partial order $\leq_K$ on $Y$, $H\subset K$  be a $K$-convex set with $0\not\in {\rm vcl}(H+K)$. Let $f:\, X\rightarrow  2^Y$ be a set-valued map, $S$ be a nonempty closed subset of $X$, $x_0\in S$ and $\gamma >0$  be a constant. For any $x\in S$, put
$$S(x) :=\,\{x^{\prime}\in S:\, f(x) \subset f(x^{\prime}) + \gamma\, d(x, x^{\prime}) H +K\}.$$ Suppose that the following conditions are satisfied:

{\rm (i)} \ $(X, d)$ is  $S(x_0)$-dynamically complete or $(X,d)$ is $(f, K)$-lower complete;

{\rm (ii)} \  $x_0\,\in\, AE(C_H, \epsilon)$, i.e., $f(x_0) \not\subset f(S) +\epsilon\,H+K;$

{\rm (iii)} \  for any  $x\in S(x_0)$, $S(x)$ is dynamically closed.

Then, there exists $\bar{x}\in S$ such that

{\rm (a)} \ $f(x_0)\,\subset\, f(\bar{x}) + \gamma d(x_0, \bar{x}) \,H +K;$

{\rm (b)} \ $\bar{x}$ is a  strict minimizer of the set-valued map $x\mapsto\, f(x) + \gamma  d(\bar{x}, x)\, H$, i.e.,
$$f(\bar{x})\,\not\subset\, f(x) +\gamma\, d(\bar{x}, x)\, H +K,\ \ \forall x\in S\backslash\{\bar{x}\};$$

{\rm  (c)} \  $d(x_0, \bar{x})\, H\,\cap\,(\epsilon/\gamma)({\rm cone}(C_H)\backslash C_H)\,\not=\,\emptyset.$}\\

{\sl Proof.} \  By condition (i) and $S$ being closed in $X$, we know that $(S, d)$ is  $S(x_0)$-dynamically complete.
Now, substituting $(S, d)$ for $(X, d)$ in Corollary 5.5, we conclude that  there exists $\bar{x}\in S$ such that (a) and (b)  hold.  Next, we show that (c) holds.
By (ii),  there exists $y_0\in f(x_0)$ such that $y_0\not\in f(S) +\epsilon\, H +K.$ Since  (a) holds, we have
$$y_0\in f(x_0)\subset f(\bar{x}) +\gamma\, d(x_0, \bar{x})\, H +K.$$
Thus,  there exists $\bar{y}\in f(\bar{x}),\ h_0\in H$ and $k_0\in K$  such that
$$y_0 = \bar{y} +\gamma\,d(x_0, \bar{x})\, h_0 + k_0.$$
If $d(x_0, \bar{x})\,h_0\,\in\, (\epsilon/\gamma) C_H\,=\,(\epsilon/\gamma)(H+K),$  then $\gamma\,  d(x_0, \bar{x})\, h_0\,\in\, \epsilon H +K,$
which leads to $$y_0\,\in\, f(\bar{x}) +\epsilon\, H +K +K\subset  f(S) +\epsilon\, H +K,$$
contradicting the assumption that $y_0\not\in f(S) + \epsilon\,H +K$.  Thus, $d(x_0, \bar{x})\, h_0\,\not\in\, (\epsilon/\gamma)\, C_H$.
On the other hand,  clearly $d(x_0, \bar{x})\, h_0\,\in\, d(x_0, \bar{x})\, H$ and $d(x_0, \bar{x})\, h_0\,\in\, {\rm cone}(C_H)$.
Thus, $d(x_0, \bar{x})\, h_0\,\in\, d(x_0, \bar{x})\, H \,\cap\, (\epsilon/\gamma)\, ({\rm cone}(C_H)\backslash C_H)$ and (c) holds.\hfill\framebox[2mm]{}\\

Obviously, Theorem 6.6 extends [42, Theorem 6.3] and [22, Theorem 5.11] to the case that $f$  is a set-valued map. Besides,  the condition that $(f(S)-f(x_0))\cap (-\epsilon\,{\rm cone}(C_H)\backslash C_H)$ is $K$-bounded (See [42, 22])  has been completely removed. Here, a set $M\subset Y$ is said to be $K$-bounded (by scalarization) iff for every $l\in K^+$, $\inf\{l(y):\, y\in M\} >-\infty$.

Similarly, we give another expression of Corollary 5.7 as follows.\\

{\bf Theorem 6.7.  } \  {\sl Let $(X, d)$ be a metric space, $Y$ be  a real vector space, $K\subset Y$ be a pointed convex cone specifying a partial order $\leq_K$ on $Y$, $H\subset K$ be a $K$-convex set with $0\not\in {\rm vcl}(H+K)$. Let $f:\,X\rightarrow  2^Y$
be a set-valued map, $S$ be a nonempty closed subset of $X$,   $x_0\in S$ and $\gamma >0$ be a constant.  Suppose that the following conditions are satisfied:

{\rm (i)} \ $(X,d)$ is $(f, K)$-lower complete (or complete);

{\rm (ii)} \  $x_0\,\in\, AE(C_H, \epsilon)$;

{\rm (iii)} \ $f$ is $K$-slm and has $K_H$-$v_0$-closed values, where $v_0\in H+K$.

Then, there exists $\bar{x}\in S$ such that

{\rm (a)} \ $f(x_0)\,\subset\, f(\bar{x}) + \gamma d(x_0, \bar{x}) \,H +K;$

{\rm (b)} \  $f(\bar{x})\,\not\subset\, f(x) +\gamma\, d(\bar{x}, x)\, H +K,\ \ \forall x\in S\backslash\{\bar{x}\};$

{\rm  (c)} \  $d(x_0, \bar{x})  \leq \epsilon/\gamma$.}\\

{\sl Proof.} \ From condition (i), i.e., $(X, d)$ is $(f, K)$-lower complete, we can easily deduce that $(S, d)$ is $(f,K)$-lower complete. Now, substituting $(S, d)$ for $(X, d)$ in Corollary 5.7, we immediately obtain the result. \hfill\framebox[2mm]{}\\

{\bf Corollary 6.8.} \  {\sl In Theorem 6.7, the result remains true if condition {\rm (iii)}  is replaced by the following condition:

{\rm (iii$^{\prime}$)} \  $f$ is $K$-slm, $f(x)\subset Y$ is weakly compact for all $x\in X$, $H$ is weakly compact and $K$ is closed, where  $Y$ is a locally convex space.}\\

{\sl Proof.} \  By Theorem 6.7, we only need to show that $f$ has $K_H-v_0$-closed values for some $v_0\in H+K$. Since $H$ is weakly compact, $\lambda\,H$ is  weakly compact.  Since  $K$ is a closed convex cone, $K$  is weakly closed. Thus, $\lambda\, H+K$ is weakly closed. Finally, $f(x)$ is weakly compact, so  $f(x) +\lambda\, H +K$ is weakly closed. Certainly, $f(x) +\lambda\,H +K$ is $v_0$-closed. \hfill\framebox[2mm]{}\\

{\bf Corollary   6.9. } \  {\sl In Theorem 6.7, the result remains true if  condition {\rm (iii)}  is replaced by the following condition:

{\rm (iii$^{\prime\prime}$)} \  $f$ is $K$-slm, $f(x)\subset Y$ is $\sigma(Y, K^+)$-sequentially compact  for all $x\in X$, $H$ is $\sigma(Y, K^+)$-countably compact  and $K$ is closed, where  $Y$ is a locally convex space.}\\

{\sl Proof.} \  By Theorem 6.7, we only need to show that  $f$ has $K_H-v_0$-closed for some $v_0\in H+K$.
Take any $v_0\in H+K$.  Let $x\in X$ and let $z\in{\rm vcl}_{v_0}(f(x) +\lambda \,H +K)$. Then, there exists a sequence $(\lambda_n)$ with all $\lambda_n \geq 0$ and $\lambda_n \rightarrow 0$ such that
$$ z+\lambda_n v_0\,\in\, f(x) +\lambda\, H +K.$$
Thus, for every $n$, there exists $y_n\in f(x)$, $h_n\in H$ and $k_n\in K$ such that
$$z+\lambda_n v_0\,=\, y_n +\lambda\, h_n + k_n.\eqno{(6.3)}$$
Since $f(x)$ is $\sigma(Y, K^+)$-sequentially compact, there exists a subsequence $(y_{n_i})_i$ of $(y_n)_n$ and $y_0\in f(x)$ such that
$$y_{n_i}\,\rightarrow\, y_0\ \  (i\rightarrow\infty)\ \ \   {\rm in} \ \ (Y, \sigma(Y, K^+)).\eqno{(6.4)}$$
Since the sequence $(h_{n_i})_i\,\subset\, H$ and $H$ is $\sigma(Y, K^+)$-countably compact, there exists $h_0\in H$  such that $h_0$ is a $\sigma(Y, K^+)$-cluster point of $(h_{n_i})_i$. For any $l\in K^+$,
$l(h_0)$  is a cluster point of the real number sequence $(l(h_{n_i}))_i$. Thus,  there exists a subsequence $(l(h_{n_{i_j}}))_j$ of the sequence $(l(h_{n_i}))_i$ such that
$$(l(h_{n_{i_j}}))_j\, \rightarrow\, l(h_0)\  \   (j\rightarrow\infty).\eqno{(6.5)}$$
By (6.4), we have
$$l(y_{n_i})\, \rightarrow\, l(y_0)\       \ (i\rightarrow\infty)\ \ {\rm and} \ \ l(y_{n_{i_j}})\,\rightarrow\, l(y_0)\ \ (j\rightarrow\infty). \eqno{(6.6)}$$
From (6.3), we have
$$l(z) +\lambda_{n_{i_j}} l(v_0)\,=\, l(y_{n_{i_j}}) +\lambda\, l(h_{n_{i_j}}) + l(k_{n_{i_j}})\,\geq\, l(y_{n_{i_j}}) +\lambda\, l(h_{n_{i_j}}).\eqno{(6.7)}$$
Letting $j\rightarrow\infty$ and combining (6.5), (6.6) and (6.7), we have
$$l(z)\,\geq\, l(y_0) +\lambda\,l(h_0) \ \  {\rm  and} \ \   l(z-y_0-\lambda\, h_0)\,\geq\,0.$$
From this, $$z- y_0 -\lambda\, h_0\,\in\, K^{++} \,=\,K.$$
That is, $$z\,\in\, y_0+\lambda\, h_0 +K\,\subset\,  f(x) +\lambda\, H +K.$$
Thus, $f(x) +\lambda\, H +K$ is $v_0$-closed. \hfill\framebox[2mm]{}\\

\section*{\large\bf  7   \   Set-valued  EVP where  perturbation contains $\sigma$-convex set  }

\hspace*{\parindent}
Vector-valued EVPs, where the objective functions are a vector-valued map $f:\, X\rightarrow Y$ and the perturbations are of type $d(\cdot, \cdot)\,H$, where $H$ is a  $\sigma$-convex set, have been considered by Bednarczuk and Zagrodny [7], Tammer and Z$\breve{a}$linescu [51] and Qiu [42]; for details, see [7, Theorem 4.1], [51, Theorem 6.2] and [42, Theorem 6.8]. Moreover, Qiu [44, Theorem 4.2]  gave such an EVP, where the objective function is a set-valued map $f:\, X\rightarrow 2^Y$ and the perturbation is as the above form, i.e., $d(\cdot, \cdot)\, H$.  There, the lower boundedness  condition is as follows: there exists $k^*\in K^+\cap H^{+s}$  such that  $k^*$ is lower bounded on $f(S(x_0))$. Obviously, it is not the weakest (See Section 3). In this section, under the weakest lower boundedness condition we shall give a set-valued EVP, where the objective function is a set-valued map and the perturbation contains a $\sigma$-convex set.

First, we recall some facts on  $\sigma$-convex sets. Let $Y$ be a t.v.s. and $B\subset Y$ be nonempty.  A convex series of points of $B$ is a series of the form $\sum_{n=1}^{\infty} \lambda_n b_n$, where every $b_n\in B$, every $\lambda_n \geq 0$  and $\sum_{n=1}^{\infty}\lambda_n\,=\,1$. $B$  is said to be a $\sigma$-convex  set  iff every convex series of its points converges to a point of $B$ (see [39, 43]).
 It is easy to show that a set is $\sigma$-convex iff it is cs-complete and bounded (concerning cs-complete sets, see [51, 53]). Suppose  that $B$ is a $\sigma$-convex set. Then, for a sequence $(b_n)$ in $B$ and a real sequence $(\lambda_n)$ with $\lambda_n\geq 0$ and $0 < \sum_{n=1}^{\infty} \lambda_n < +\infty$,
$\sum_{n=1}^{\infty} \lambda_n b_n/\sum_{n=1}^{\infty}\lambda_n$ is a convex series in $B$ and it converges to some point $\bar{b}\in B$.
Thus, $\sum_{n=1}^{\infty}\lambda_n b_n$ converges to $(\sum_{n=1}^{\infty}\lambda_n) \bar{b}\,\in\,(\sum_{n=1}^{\infty}\lambda_n) B$.
 A set $B$ in $Y$ is said to be sequentially complete if every Cauchy sequence  $(b_n)$ in $B$, converges to a point of $B$.
 In [7], a sequentially complete set  is called a semi-complete set. It is easy to show that every sequentially complete, bounded  convex set  is a $\sigma$-convex set (see  [51, Remark 6.1]).
If $Y$ is a locally complete locally convex space (see [39, 49]),  then every locally closed, bounded convex set in $Y$ is a $\sigma$-convex set.  However, a $\sigma$-convex set needn't be sequentially complete or even needn't be sequentially closed (or needn't be locally closed). In fact, an open ball in a Banach space is $\sigma$-convex, but it isn't closed. For details, see [39, 43].

The following theorem extends and improves Bednarczuk and Zagrodny [7, Theorem 4.1], Tammer and Z$\breve{a}$linescu [51, Theorem 6.2], Liu and Ng [35, Theorem 3.5(ii)],  Guti\'{e}rrez, Jim\'{e}nez and Novo [22, Theorem 5.12], Qiu [42, Theorem 6.8] and Qiu [44, Theorem 4.2].\\

{\bf Theorem 7.1.} \  {\sl Let $(X, d)$ be a metric space, $x_0\in X$, $Y$ be a t.v.s. quasi-ordered by a convex cone $K$, $H\subset K$ be a $\sigma$-convex set such that $0\not\in {\rm vcl}(K+H)$ and let $f:\, X\rightarrow 2^Y$  be a set-valued map.
For any $x\in X$, put $$S(x):=\{x^{\prime}\in X:\, f(x) \subset f(x^{\prime}) +d(x, x^{\prime})\, H +K\}.$$

Suppose that the following conditions are satisfied:

{\rm (i)} \ $(X, d)$ is $S(x_0)$-dynamically complete;

{\rm (ii)} \  there exists $\epsilon >0$ such that $f(x_0)\, \not\subset\, f(S(x_0)) +\epsilon\,H +K;$

{\rm (iii)} \ $f$ is $K$-slm and has $K$-closed values, i.e., for any $x\in X$, $f(x) + K$ is closed.

Then, there exists $\bar{x}\in X$ such that

{\rm (a)} \ $f(x_0)\,\subset\, f(\bar{x}) + d(x_0, \bar{x})\, H +K;$

{\rm (b)} \  $\forall x\in X\backslash\{\bar{x}\},\ f(\bar{x})\,\not\subset\, f(x) + d(\bar{x}, x)\, H +K.$}\\

{\sl Proof.} \  By Theorem 5.3,  we only need to prove  that for any $x\in S(x_0)$, $S(x)$ is dynamically closed.
Let $(x_n) \subset  S(x) \subset S(x_0)$, $x_{n+1} \in S(x_n),\ \forall n,$ and $x_n \rightarrow u$. Take any $n_0\in {\bf N}$ and put $z_1:= x_{n_0}$. We assume that $z_1\not= x_0$.
We may take a subsequence $(z_n)$  from $(x_k)$ such that $z_{n+1}\in  S(z_n)$ and $z_{n+1} \not= z_n,\ \forall n.$ Obviously, we have $d(z_{n+1}, u) \rightarrow 0\ (n\rightarrow\infty)$.
By condition (ii), there exists $y_0 \in f(x_0)$ such that
$$y_0\,\not\in\, f(S(x_0)) +\epsilon\, H +K.$$
Since $z_1 = x_{n_0}\,\in\, S(x_0)$, we have
$$y_0\,\in\, f(x_0)\,\subset\, f(z_1) + d(x_0, z_1)\, H +K.$$
Thus, there exists $y_1\in f(z_1)$ such that
$$y_0\,\in\, y_1 + d(x_0, z_1)\,H +K.\eqno{(7.1)}$$
Since  $z_2\,\in\, S(z_1)$, we have
$$y_1\,\in\, f(z_1)\,\subset\, f(z_2) + d(z_1, z_2)\, H +K.$$
Thus, there exists $y_2\,\in\, f(z_2)$ such that
$$y_1\,\in\, y_2 + d(z_1, z_2)\, H +K. \eqno{(7.2)}$$
In general, Let $y_n\,\in\, f(z_n)$ be given.
Since $z_{n+1} \in\, S(z_n)$, we have
$$y_n\,\in\, f(z_n)\,\subset\, f(z_{n+1}) + d(z_n, z_{n+1})\, H +K.$$
Thus, there exists $y_{n+1}\,\in\, f(z_{n+1})$ such that
$$y_n\,\in\, y_{n+1} + d(z_n, z_{n+1})\, H +K. \eqno{(7.3)}$$
By (7.1), (7.2) and (7.3), we have
$$y_0 +y_1 +\cdots + y_n\,\in\, y_1 + y_2 +\cdots +y_{n+1} + \left(d(x_0, z_1)+\sum_{i=1}^n d(z_i, z_{i+1})\right)\,H +K.$$
Thus,
$$y_0 \,\in\,  y_{n+1} + \left(d(x_0, z_1)+\sum_{i=1}^n d(z_i, z_{i+1})\right)\,H +K.$$
From this,
$$y_{n+1} - y_0\,\in\, -\left(d(x_0, z_1) + \sum_{i=1}^n d(z_i, z_{i+1})\right)\,H -K$$
and
$$\xi_H(y_{n+1} - y_0)\,\leq\, -d(x_0, z_1) - \sum_{i=1}^n d(z_i, z_{i+1}).\eqno{(7.4)}$$
Since  $y_{n+1}\,\in\, f(z_{n+1})\,\subset\, f(S(x_0))$ and $y_0\,\not\in\, f(S(x_0)) + \epsilon\,H +K,$
we have
$$y_0\,\not\in\, y_{n+1} +\epsilon\, H +K\ \ \ \ \  {\rm  and}\ \ \ \ \
y_{n+1} -y_0\,\not\in\, -\epsilon\, H -K.$$
Thus,
$$\xi_H(y_{n+1} -y_0)\,\geq \,-\epsilon.\eqno{(7.5)}$$
Combining (7.4) and (7.5), we have
$$-\epsilon\,\leq\, -d(x_0, z_1) - \sum_{i=1}^n d(z_i, z_{i+1}).$$
Thus, $$d(x_0, z_1) + \sum_{i=1}^n d(z_i, z_{i+1})\,\leq\, \epsilon$$ and $$\sum_{i=1}^{\infty} d(z_i, z_{i+1})\, <\, +\infty.\eqno{(7.6)}$$
Now, for any $y_1^{\prime} \in f(z_1)$, we  have
$$y_1^{\prime}\,\in\, f(z_1) \,\subset\, f(z_2)  + d(z_1, z_2)\, H +K.$$
So there exists $y_2^{\prime}\,\in\, f(z_2),\ h_1\in H$ such that
$$y_1^{\prime}\,\in\, y_2^{\prime} + d(z_1, z_2)\, h_1 +K.\eqno{(7.7)}$$
For $y_2^{\prime}\,\in\, f(z_2) \,\subset\, f(z_3) + d(z_2, z_3)\, H +K, $
there exists $y_3^{\prime}\in f(z_3),\ h_2\in H$ such that
$$y_2^{\prime}\,\in \, y_3^{\prime} + d(z_2, z_3)\, h_2 +K.\eqno{(7.8)}$$
In general, if $y_n^{\prime} \in f(z_n)$ is given, then there exists $y_{n+1}^{\prime} \in f(z_{n+1}), \ h_n\in H$ such that
$$y_n^{\prime}\,\in\, y_{n+1}^{\prime} + d(z_n, z_{n+1})\,h_n  +K.\eqno{(7.9)}$$
By (7.7), (7.8) and (7.9), we have
$$y_1^{\prime} + y_2^{\prime} +\cdots  y_n^{\prime}\,\in\, y_2^{\prime} + y_3^{\prime} +\cdots  +y_{n+1}^{\prime} +\sum_{i=1}^n d(z_i, z_{i+1})\, h_i+K$$
From this,
$$y_1^{\prime}\,\in\, y_{n+1}^{\prime} +\sum_{i=1}^n d(z_i, z_{i+1}) h_i +K\,=\, y_{n+1}^{\prime} +\left(\sum_{i=1}^n d(z_i, z_{i+1})\right) h_n^{\prime} +K,\eqno{(7.10)}$$
where $h_n^{\prime}\,=\,\sum_{i=1}^n d(z_i, z_{i+1}) h_i/\sum_{j=1}^n d(z_j, z_{j+1})\,\in\, H.$
Since $\sum_{i=1}^n d(z_i, z_{i+1})\,\geq\, d(z_1, u) - d(z_{n+1}, u)$,  by (7.10) and using the condition $f$ being $K$-slm, we have
\begin{eqnarray*}
y_1^{\prime}\,&\in&\, y_{n+1}^{\prime} + \left(\sum_{i=1}^n d(z_i, z_{i+1})\right)h_n^{\prime}+K\\
&\subset&\,y_{n+1}^{\prime} + (d(z_1, u) - d(z_{n+1}, u)) h_n^{\prime} +K\\
  &\subset&\,f(z_{n+1})+ (d(z_1, u) -d(z_{n+1}, u)) h_n^{\prime} + K\\
&\subset&\, f(u) + (d(z_1, u) - d(z_{n+1}, u)) h_n^{\prime} + K\\
&=&\, f(u) + d(z_1, u)h_n^{\prime} - d(z_{n+1}, u) h_n^{\prime} + K.\hspace{6cm} (7.11)
\end{eqnarray*}
Since $H$ is $\sigma$-convex, by (7.6) we have
$$ h_n^{\prime}\, \rightarrow\, \bar{h}:=\sum_{i=1}^{\infty} d(z_i, z_{i+1})\, h_i/\sum_{j=1}^{\infty} d(z_j, z_{j+1})\,\in \, H.$$
Since $(h_n^{\prime})$ is bounded and $d(z_{n+1}, u) \rightarrow 0$, we have $d(z_{n+1}, u) h_n^{\prime}\,\rightarrow\, 0\  (n\rightarrow\infty)$.
Now, by (7.11) and remarking that $f(u) +K$ is closed, we have
$$y_1^{\prime} \,\in\, f(u) + d(z_1, u)\,\bar{h} + K\,\subset\, f(u) + d(z_1, u)\, H +K.$$
This holds for all $y_1^{\prime} \in f(z_1)$. Thus,
$$f(z_1)\,\subset\, f(u) + d(z_1, u)\,H +K.$$
That is,
$u\,\in S(z_1) = S(x_{n_0})\,\subset\, S(x).$\hfill\framebox[2mm]{}\\

Assume that $Y$ is a locally convex space. Then, from the above proof, we see that $(d(z_{n+1}, u)\, h^{\prime}_n)_n$  is locally convergent to $0$ and $(h^{\prime}_n)_n$ is locally convergent to $\bar{h}$. Thus, from (7.11) we know that the assumption that $f$ has $K$-closed values (See   condition (iii)) can be replaced by a slightly weaker  one: $f$ has $K$-locally closed values, i.e.,  for any $x\in X$, $f(x) +K$ is locally closed.\\

\noindent {\bf Acknowledgements}

 This work was supported by the National Natural Science Foundation of China (Grant Nos. 11471236, 11561049).\\

\noindent{\bf References} \vskip 10pt
\begin{description}
\footnotesize
\item{[1]}  M. Adan,   V. Novo, Weak efficiency in vector optimization using a closure of algebraic type under cone-convexlikeness, European J. Oper. Res., 149 (2003), pp. 641-653.

\item{[2]} M. Adan,   V. Novo,  Proper efficiency in vector optimization
on real linear spaces, J. Optim. Theory Appl., 121 (2004), pp. 515-540.

\item{[3]}  Y. Araya,  Ekeland's variational principle and its
equivalent theorems in vector optimization, J. Math. Anal. Appl.,
 346 (2008), pp. 9-16.

\item{[4]} T. Q. Bao, B. S. Mordukhovich, Variational principles for
set-valued mappings with applications to multiobjective
optimization, Control Cybern., 36 (2007), pp. 531-562.

\item{[5]} T. Q. Bao, B. S. Mordukhovich, Relative Pareto minimizers
for multiobjective problems: existence and optimality conditions,
Math. Program, Ser.A, 122 (2010), pp. 301-347.

\item{[6]}  E. M. Bednarczuk,  M. J. Przybyla, The vector-valued
variational principle in Banach spaces ordered by cones with
nonempty interiors, SIAM J. Optim. 18 (2007), pp. 907-913.

\item{[7]}  E. M. Bednarczk, D. Zagrodny, Vector variational
principle, Arch. Math. (Basel), 93 (2009), pp. 577-586.

\item{[8]}  G. Y. Chen, X. X. Huang, X. G. Yang, Vector Optimization, Set-Valued
and Variational Analysis, Springer-Verlag, Berlin, 2005.

\item{[9]} Y. Chen, Y. J. Cho, L. Yang, Note on the results with  lower semi-continuity, Bull. Korean Math. Soc., 39 (2002), pp. 535-541.

\item{[10]} D. Dentcheva, S. Helbig, On variational principles,
level sets, well-posedness, and $\epsilon$-solutions in vector
optimization, J. Optim. Theory Appl., 89 (1996), pp. 325-349.

\item{[11]}  W. S. Du, On some nonlinear problems induced by an abstract maximal element principle, J. Math. Anal. Appl., 347 (2008), pp. 391-399.

\item{[12]}  I. Ekeland, Sur les prob\`{e}mes variationnels, C. R. Acad. Sci.
Paris,  275 (1972), pp. 1057-1059.

\item{[13]}  I. Ekeland,  On the variational principle, J. Math.
Anal. Appl.,  47 (1974), pp. 324-353.

\item{[14]} I. Ekeland, Nonconvex minimization problems, Bull. Amer.
Math. Soc., (N.S.) 1 (1979), pp. 443-474.

\item{[15]}  C. Finet, L. Quarta, C. Troestler,  Vector-valued
variational principles, Nonlinear Anal., 52 (2003), pp. 197-218.

\item{[16]}  F. Flores-Baz\'{a}n, C. Guti\'{e}rrez,  V. Novo,  A
Br\'{e}zis-Browder principle on partially ordered spaces and related
ordering theorems, J. Math. Anal. Appl., 375 (2011), pp. 245-260.

\item{[17]} Chr. Gerstewitz (Tammer), Nichtkonvexe Dualit\"{a}t in
der Vektoroptimierung, Wiss. Z. TH Leuna-Merseburg 25 (1983), pp.
357-364.

\item {[18]}  A. G\"{o}pfert, H. Riahi,  Chr. Tammer,  C.
Z$\breve{a}$linescu,  Variational Methods in Partially Ordered
Spaces, Springer-Verlag, New York, 2003.

 \item{[19]} A. G\"{o}pfert,  C. Tammer,  C. Z$\breve{a}$linescu,  On the
vectorial Ekeland's variational principle and minimal point theorems
in product spaces, Nonlinear Anal.,  39 (2000), pp.  909-922.

\item{[20]} C. Guti\'{e}rrez, B. Jim\'{e}nez, V. Novo, On
approximate efficiency in multiobjective programming, Math. Methods
Oper. Res., 64(2006),  pp. 165-185.

\item{[21]} C. Guti\'{e}rrez, B. Jim\'{e}nez, V. Novo,  A unified
approach and optimality conditions for approximate solutions of
vector optimization problems, SIAM J. Optim., 17(2006),  pp. 688-710.

\item{[22]} C. Guti\'{e}rrez, B. Jim\'{e}nez, V. Novo, A set-valued
Ekeland's variational principle in vector  optimization, SIAM J.
Control. Optim., 47 (2008), pp.  883-903.

\item{[23]} C. Guti\'{e}rrez,  V. Novo, J. L. R\'{o}denas-Pedregosa,  T. Tanaka, Nonconvex separation functional in linear spaces with  applications to  vector equilibria, SIAM J. Optim., 26 (2016), pp. 2677-2695.

\item{[24]} T. X. D. Ha, Some variants of the Ekeland  variational
principle for a set-valued map, J. Optim. Theory Appl., 124 (2005),
pp. 187-206.

\item{[25]} A. H. Hamel, Equivalents to Ekeland's variational principle in uniform spaces, Nonlinear Anal., 62 (2005), pp. 913-924.

\item{[26]} F. He, J. H. Qiu, Sequentially lower complete spaces and Ekeland's variational principle, Acta Math. Sin. (Engl. Ser.),  31 (2015), pp. 1289-1302.

\item{[27]} R. B. Holmes, Geometric Functional Analysis and its Applications, Springer, New York, 1975.

\item{[28]} G. Isac, The Ekeland's principle and the
Pareto $\epsilon$-efficiency, In: M. Tamiz (ed.) Multi-Objective
Programming and Goal Programming: Theories and Applications, Lecture
Notes in Econom. and Math. Systems, vol. 432, Springer-Verlag,
Berlin, 1996, pp. 148-163.

\item{[29]} P. Q. Khanh,  D. N. Quy, Versions of Ekeland's variational principle involving set perturbations, J. Glob. Optim.,  57 (2013), pp. 951-968.

\item{[30]}  J. L. Kelley, I. Namioka, W. F. Donoghue,Jr., K. R. Lucas, B. J. Pettis, E. T.  Poulsen, G. B. Price,  W. Robertson, W. R. Scott, K. T. Smith,  Linear Topological
Spaces, Van\ Nostrand, Princeton, 1963.

\item{[31]}  G.  K\"{o}the, , Topological Vector Spaces I, Springer-Verlag,
Berlin, 1969.

\item{[32]}  D. Kuroiwa, On set-valued optimization, Nonlinear Anal., 47 (2001), pp. 1395-1400.

\item{[33]} S. J. Li, X. Q. Yang,  G. Y.  Chen, Vector Ekeland  variational principle, in: Vector Variational Inequalities and Vector Equilibria,  F. Giannessi (ed.), Nonconvex Optimization and its Applications, vol. 38, Klower, Dordrecht, 2000, pp. 321-333.

\item{[34]} S. J. Li, W. Y. Zhang, A minimization theorem for a set-valued mapping, Applied Math. Letters, 21 (2008), pp. 769-773.

\item{[35]}  C. G. Liu, K. F. Ng, Ekeland's variational principle for
set-valued functions, SIAM J. Optim., 21 (2011), pp. 41-56.

\item{[36]} D. T. Luc, Theory of Vector Optimization, Lecture Notes in Econom. and Math. System 319, Springer, Berlin, 1989.

\item{[37]} D. G. Luenberger, New optimality principles for economic efficiency and equilibrium, J. Optim. Theory  Appl., 75 (1992), pp. 221-264.

\item{[38]}  A. B. N\'{e}meth, A nonconvex vector minimization problem,
Nonlinear Anal.,  10 (1986), pp.  669-678.

\item{[39]}  P. P\'{e}rez Carreras, J. Bonet, Barrelled Locally Convex Spaces,
North-Holland, Amsterdam, 1987.

\item{[40]}  J. H. Qiu, A generalized Ekeland vector variational
principle and its applications in optimization, Nonlinear Anal., 71
(2009), pp. 4705-4717.

\item{[41]} J. H. Qiu, On Ha's version of set-valued Ekeland's
variational principle, Acta Math. Sin. (Engl. Ser.), 28 (2012),
pp. 717-726.

\item{[42]} J. H. Qiu, Set-valued quasi-metrics and a general
Ekeland's variational principle in vector optimization, SIAM J.
Control  Optim., 51 (2013), pp. 1350-1371.

\item{[43]} J. H. Qiu, The domination property for efficiency and
Bishop-Phelps theorem in locally convex spaces, J. Math. Anal.
Appl., 402 (2013), pp. 133-146.

\item{[44]} J. H. Qiu, A pre-order principle and set-valued  Ekeland variational principle, J. Math. Anal. Appl., 419 (2014), pp. 904-937.

\item{[45]} J. H. Qiu, An equilibrium version of vectorial Ekeland variational principle and its applications to equilibrium problems, Nonlinear Analysis: RWA, 27 (2016), pp. 26-42.

\item{[46]} J. H. Qiu, An equilibrium version of set-valued Ekeland variational principle and its applications to set-valued vector equilibrium problems, Acta Math. Sin. (Engl. Ser.), 33 (2017), pp. 210-234.

\item{[47]} J. H. Qiu, F. He, A  general vectorial Ekeland's
variational principle with a p-distance, Acta Math. Sin.  (Engl. Ser.),  29
(2013), pp. 1655-1678.

\item{[48]} J. H. Qiu, B. Li, F. He, Vectorial Ekeland's variational
principle with a w-distance and its equivalent theorems, Acta Math.
Sci. Ser. B,  32 (2012), pp. 2221-2236.

\item{[49]} S. A. Saxon, L. M. S\'{a}nchez Ruiz, Dual local
completeness, Proc. Amer. Math. Soc., 125 (1997), pp. 1063-1070.

\item{[50]} C. Tammer, A generalization of Ekeland's variational principle,
Optimization,  25 (1992), pp.  129-141.

\item{[51]} C. Tammer, C. Z$\breve{a}$linescu, Vector variational
principle for set-valued functions, Optimization, 60 (2011), pp.
839-857.

\item{[52]}  A. Wilansky, Modern Methods in Topological Vector Spaces,
McGraw-Hill, New York, 1978.

\item{[53]} C. Z$\breve{a}$linescu, Convex Analysis  in   General
Vector Spaces, World Sci., Singapore, 2002.

\item{[54]} J. Zhu, C. K. Zhong, Y. J. Cho, A generalized variational principle
and vector optimization, J. Optim. Theory Appl., 106 (2000), pp. 201-218.

\end{description}
\end{document}